\documentclass[12pt,twoside]{article}
\usepackage{amssymb}
\usepackage{color}
\definecolor{Old}{gray}{.50}
\definecolor{Grey}{gray}{.80}
\usepackage{amsmath}
\usepackage{latexsym}
\usepackage{epsfig}

\newtheorem{theorem}{Theorem}[section]
\newtheorem{definition}[theorem]{Definition}

\newtheorem{lemma}[theorem]{Lemma}
\newtheorem{proposition}[theorem]{Proposition}

\renewcommand{\thetheorem}{\thesection.\arabic{theorem}}

\renewcommand{\theequation}{\thesection.\arabic{equation}}

\newcommand{\IR}{\mathbb{R}}
\newcommand{\IZ}{\mathbb{Z}}
\newcommand{\IP}{\mathbb{P}}
\newcommand{\IE}{\mathbb{E}}
\newcommand{\IN}{\mathbb{N}}
\newcommand{\Ai}{\text{Ai}}
\newcommand{\dd}{\text{d}}

\newcommand{\tr}{\text{tr}}
\newcommand{\id}{\text{1\hspace{-0.25em}l\,}}
\renewcommand{\ss}{\text{s}}
\def\Ord{\mathcal{O}}
\begin{document}
\begin{titlepage}
\begin{center}
\vspace*{2,5cm}
 {\Large\bf  Scale Invariance of the PNG Droplet and the Airy\\[3mm] Process}
\bigskip\bigskip\bigskip\\
{\large Michael Pr\"{a}hofer and Herbert Spohn}\bigskip\\
Zentrum Mathematik and Physik Department, TU M\"unchen,\\
D-80290 M\"unchen, Germany\medskip\bigskip\\
emails: {\tt praehofer@ma.tum.de}, {\tt spohn@ma.tum.de}\\[5mm]
{\it Dedicated with admiration to David Ruelle and Yasha Sinai at the
 occasion\\ of their $65$-th birthday.}
\end{center}

\baselineskip=24pt \vspace{3cm}
\begin{center}
  {\bf Abstract}
\end{center}
We establish that the static height fluctuations of a particular growth
model, the PNG droplet, converges upon proper rescaling to a limit
process, which we call the Airy process $A(y)$.
The Airy process is stationary, it has continuous
sample paths, its single ``time'' (fixed $y$) distribution is the
Tracy-Widom distribution of the largest eigenvalue of a GUE random
matrix, and the Airy process has a slow decay of correlations as
$y^{-2}$. Roughly the Airy process describes the last line of
Dyson's Brownian motion model for random matrices. Our
construction uses a multi-layer version of the PNG model, which
can be analyzed through fermionic techniques.
Specializing our result to a fixed value of $y$, one reobtains
the celebrated result of Baik,
Deift, and Johansson on the length of the longest increasing
subsequence of a random permutation.

\vspace{1cm}

\end{titlepage}

\section{The PNG droplet} \label{sec.a}
\setcounter{equation}{0}

The polynuclear growth (PNG) model is a simplified model for layer
by layer growth \cite{Mea, Pra}. Initially one has a perfectly
flat crystal in contact with its supersaturated vapor. Once in a
while a supercritical seed is formed, which then spreads
laterally by further attachment of particles at its perimeter
sites. Such islands coalesce if they are in the same layer and
further islands may be nucleated upon already existing ones. The
PNG model ignores the lateral lattice structure and assumes that
the islands are circular and spread at constant speed.

In this paper we study the one-dimensional version in the
particular geometry where nucleation only above the ground layer
$[-t, t]$ is allowed. To be precise: at time $t$ the height is given
by the (random) function $h (x,t)$, $x \in \IR$. One requires $h(x,t)
= 0$ for $|x| > t$, in particular $h(x,0)=0$ for $x\in\IR$. At a given
time $t$, $h(x,t)$ is 
piecewise constant, takes nonnegative integer values, and has
jumps of size $\pm 1$ only. The jumps are called steps and we
distinguish between up-steps (jump size $1$) 
and down-steps (jump size $-1$). The dynamics has a
deterministic piece, according to which down-steps move with
velocity $1$ and up-steps with velocity $-1$. Surface steps disappear
upon collision. In addition there are nucleation events by which new 
steps are created. Randomly in time, $h(x,t)$ is changed to
the new profile $\tilde{h}(x,t)$ such that for the increment $\delta
h(x,t)=\tilde{h}(x,t) - h(x,t)$ one has $ \delta h(x,t) = 1$ at some
random point $x'$, $|x'| \leq t$, and 
$\delta h (x,t) = 0$ otherwise. Immediately after this nucleation
event the deterministic evolution is followed until the next
nucleation. Growing for a while the typical height 
profile has the shape of a droplet, $h(x,t) \simeq 2 \sqrt{t^2 - x^2}$
for $|x|\leq t$. Our interest are the statistical properties of
the deviations from this average shape.\\ 
{\it Warning}: For the PNG
model one has to specify the step speed and the intensity of the
nucleation events. They can be adjusted to an arbitrary value
through a linear scale change of space-time $(x,t)$. Geometrically
velocity one is distinguished and therefore adopted here. Intensity
one for nucleation events seems also natural, but in fact introduces a
string of factors 
$\sqrt 2$. Therefore we deviate from previous conventions and set
the intensity to be equal to $2$.

The one-dimensional PNG model is just one model within the KPZ
universality class for growth. However, for this model we have
very refined statistical information, the most surprising
breakthrough being the result of Baik, Deift, and Johansson
\cite{Bai}, which states that
\begin{equation}\label{a.a}
\lim_{t \rightarrow \infty} t^{-1/3}(h(0,t) - 2t) = \chi_2
\end{equation}
in distribution, cf.~for the connection to the PNG model at the end of
the introduction. $\chi_2$ has the same distribution as the largest
eigenvalue of a $N \times N$ random hermitian matrix (GUE) in the
limit of $N \rightarrow \infty$. As discovered by Tracy and Widom
\cite{Tra} the distribution function $F_2(x) = \IP (\chi_2 \leq
x)$ is governed by the Painlev\'e II equation. One has $F_2(x) =
e^{-g(x)}$, where $g'' = u^2$, $g(x)\to0$ as $x\to\infty$, and $u$ is
the global positive solution of 
$u'' = 2u^3 + xu$ (Painlev\'e II). Its asymptotics are $u(x) \simeq
\sqrt{- x/2}$ for $x\to-\infty$, $u(x)\sim\Ai(x)$ for $x\to\infty$
with $\Ai$ the Airy function. In fact through a 
simple linear transformation \cite{Pra} one obtains the height
fluctuations for any $y$, $|y|<1$, as
\begin{equation}\label{a.b}
\lim_{t \rightarrow \infty} t^{-1/3}\big(h(yt, t) - 2t \sqrt{1 - y^2}\big)
= (1- y^2)^{1/3} \chi_2 \,.
\end{equation}

At the next level of precision, one might wonder about joint
distributions of $h(x,t)$ at the same time, say of $\{h(x_1, t), h
(x_2, t) \}$. If $x_1 = y_1 t$ and $x_2 = y_2 t$, $y_1 \neq y_2$,
then the heights are so far apart that they become statistically
independent. Thus we have to consider closer by reference points.
From the KPZ theory one knows that the lateral fluctuations live
on the scale $t^{2/3}$. Therefore the natural object is
\begin{equation}\label{a.c}
y \mapsto t^{-1/3}\big(h (yt^{2/3},t)-2t\big) = h_t (y)
\end{equation}
considered as a stochastic process in $y$. For large $t$ we have
$\langle h_t(y)\rangle \simeq -y^2$ and if in (\ref{a.b}) we replace $yt$ by
$yt^{2/3}$ we obtain
\begin{equation}\label{a.d}
\lim_{t \rightarrow \infty} h_t (y) = - y^2 + \chi_2
\end{equation}
in distribution, which suggests that $h_t(y) + y^2$ tends to a
stationary stochastic process. As our main result we establish that
this is indeed the case and rather explicitly identify the
limit process. For reasons which will become clear in the sequel we
call the novel limit process the Airy process and denote it by
$A(y)$. We refer to Section \ref{sec.d} for its definition and to
Section \ref{sec.e} 
for some of its properties. Somewhat compressed, one considers
independent fermions in one space dimension as governed by the
one-particle Hamiltonian
\begin{equation}\label{a.e}
H = - \frac{d^2}{du^2} + u \,.
\end{equation}
The fermions are in their ground state at zero chemical potential,
which is the quasifree state determined by the spectral projection
onto $\{H\leq0\}$. Because of the linearly increasing potential there
is a last fermion, which has the Tracy-Widom $\chi_2$ as positional
distribution. Extending to the Euclidean space-time through
the propagator $e^{-yH}$, the fermions move along non-intersecting
world lines with some suitable statistical weight. The Airy process
$A(y)$ is the position of the last fermion at fermionic time $y$.
\begin{theorem}
\label{the0}
Let $A(y)$ be the stationary Airy
process. Then in the sense of weak convergence of finite
dimensional distributions
\begin{equation}\label{a.f}
\lim_{t \rightarrow \infty} h_t(y) = A(y) - y^2\,.
\end{equation}
\vspace{0mm}
\end{theorem}
We recover the result in \cite{Bai} as the special case of the
convergence of the distributions for some fixed value of $y$. Even
for fixed $y$ we provide an alternative proof based on one-dimensional
fermionic field theories. 

As shown in \cite{Pra, Pra2} the PNG droplet is isomorphic to a
directed polymer with Poissonian point potential which in turn is
related to the length of the longest increasing subsequence in a
random permutation. For the sake of completeness we repeat our
result in this language. We consider Poisson points of intensity 2
in the positive quadrant. A directed polymer is a piecewise linear
path, $\omega$, from $(0,0)$ to $(t + x, t-x)$, $|x| \leq t$,
with each segment of $\omega$ bordered by Poisson points, under
the constraint that their slope is in $[0,\infty]$. The
length, $L(x, t, \omega)$, of the directed polymer $\omega$ is the
number of Poisson points visited by $\omega$. We set
\begin{equation}\label{a.g}
L(x,t) = \max_{\omega} L(x,t, \omega)\,,
\end{equation}
where the maximum is taken over all allowed directed paths at
a fixed configuration of Poisson points and at specified
endpoints. Then, in distribution,
\begin{equation}\label{a.h}
L(x,t) = h (x,t)\,.
\end{equation}
Therefore Theorem \ref{the0} yields the statistical properties of the
length of the optimal path in dependence on its endpoint at
transverse distance $y t^{2/3}$ away from $(t,t)$.

To give a brief outline: In the following section we introduce the
multi-layer PNG model, whose last layer is the PNG droplet. The
multi-layer PNG has the remarkable property that the distribution
at time $t$ is the uniform distribution on all admissible height
lines. Such kind of ensemble can be analyzed through Euclidean
Fermi fields. Our case maps onto independent fermions on the
lattice $\IZ$ with the usual nearest neighbor hopping energy and
subject to a linearly increasing external potential. The PNG
droplet corresponds to the last fermionic world line. The convergence
of the moments of the multi-layer PNG model in essence reduces to the
convergence of the discrete Fermi propagator to the continuum
Fermi propagator corresponding to the Hamiltonian (\ref{a.e}). In
the final section we establish some properties of the Airy process
and discuss the two-point function $\langle(h(y t^{2/3}, t) -
h(0,t))^2\rangle$. 

\section{The multi-layer PNG model} \label{sec.b}
\setcounter{equation}{0}

Inspired by the beautiful work
of Johansson on the Aztec diamond \cite{Joh} we enlarge the PNG
droplet to the multi-layer PNG model. While this looks like a further
complication, in fact the construction will provide us with a powerful
machinery to analyze the statistics of the PNG droplet.

\begin{figure}[tbp]
  \begin{center}
    \mbox{\epsfxsize=14cm\epsffile{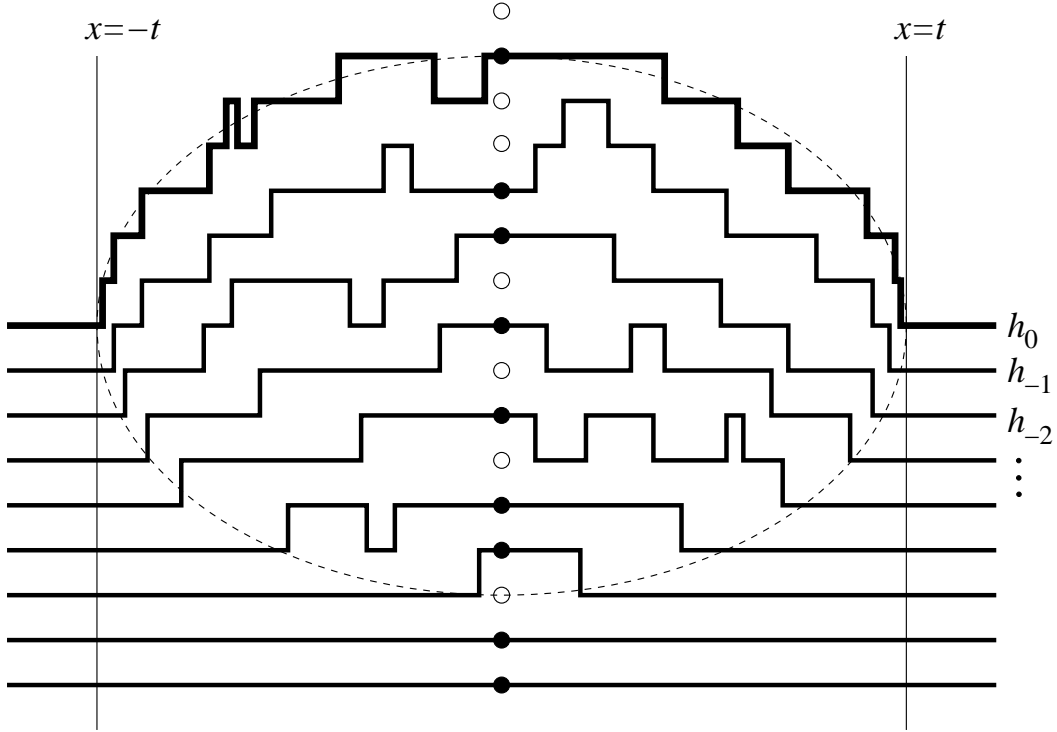}}
    \caption{\it A snapshot of a multi-layer PNG configuration at time
      $t$. As a guide to the eye the asymptotic droplet shape is
      indicated.}
    \label{fig:multiPNG}
  \end{center}
\end{figure}

Instead of a single PNG line we now consider a collection of such lines,
denoted by  $h_\ell(x,t)$, $\ell\in\IN_-=\{0,-1,-2,\dots\}$, taking
values in $\IZ$, which are arranged in ascending order,
$h_\ell(x,t)<h_{\ell+1}(x,t)$, for all $x\in\IR$, $t\geq0$, cf. Figure
\ref{fig:multiPNG} with a typical snapshot of the multi-layer PNG
model at time $t$. The lines have jump size $\pm1$ and the total
number of up- and down-steps is assumed to be finite. For each $t\geq
0$ one has $h_\ell(x,t)=\ell$ for $|x|\geq t$. To be definite we
require
$h_\ell(x,t)-\lim_{\epsilon\searrow0}h_\ell(x\pm\epsilon,t)\in\{0,1\}$,
which means that the height lines are upper semi-continuous. 
The set of all such height line configurations is denoted by
$\Lambda_t$. The height lines evolve in time under a rule which
is based on Viennot's geometric construction to prove the  
Robinson-Schensted(-Knuth) correspondence \cite{Viennot}, hence called RSK
dynamics by us. The initial line configuration is $h_\ell(x,0)=\ell$,
$\ell\in\IN_-$, for all $x$. Under the RSK dynamics only the region
$|x|< t$ is modified. The top line evolves stochastically 
like the PNG droplet. In the lower lying lines the steps move and
coalesce according to the PNG rules. However nucleations are
determined by the annihilation events in the neighboring line
above. Thus at time $t$ if in the $\ell$-th height line a collision of
an up-step and a down-step occurs at position $x$, they disappear in
this line only to reappear as nucleation at $(x,t)$ in line
$\ell-1$. Clearly the RSK dynamics respects the ordering
$h_\ell<h_{\ell+1}$.  

To describe the configuration space of the multi-layer PNG model at
time $t\geq0$, we introduce the step positions as coordinates. Let
$n_\ell$ be the number of up-steps in height line $\ell$. Since
$h_\ell(-t,t)=h_\ell(t,t)$, it must equal the number of down-steps in
height line $\ell$. If $n_\ell=0$, $h_\ell(x,t)=\ell$ for all
$x\in\IR$. If $n_\ell>0$, the position of the $j$-th up-step in height
line $\ell$ is denoted by $y^{\ell,+}_j$, 
$-t<y^{\ell,+}_1<\dots<y^{\ell,+}_{n_\ell}<t$ and the position of the  
$j$-th down-step in the same height line is denoted by $y^{\ell,-}_j$,
$-t<y^{\ell,-}_1<\dots<y^{\ell,-}_{n_\ell}<t$. 
We set ${\bf n}=(n_0,n_{-1},\dots)$, and 
$|{\bf n}|=\sum_{\ell\leq0}n_\ell$. For the RSK dynamics the
total number of steps, $2|{\bf n}|$, is finite with probability one,
since $|{\bf n}|$ equals the total number of nucleation events
up to time $t$, i.e.~the number of Poisson points in the triangle
$\{(x',t')\big|\,|x'|<t',\,t'\leq t\}$.
We denote by $\Gamma_t({\bf n})$ the set of all step configurations
$\big((y^{\ell,+}_j,y^{\ell,-}_j)_{1\leq j\leq
  n_\ell}\big)_{\ell\leq0}$ resulting from an admissible line 
configuration
$\big(h_\ell(x,t)\big)_{\ell\leq0}\in\Lambda_t$. $\Gamma_t(0)$ is a 
single point and $\Gamma_t({\bf n})$ is naturally embedded in
$[-t,t]^{2|{\bf n}|}$. By definition of $\Lambda_t$, if an
up-step and a down-step in the same height line are at the
same location, they represent a nucleation and not a collision. 
Finally $\Gamma_t=\bigcup_{|{\bf n}|<\infty}\Gamma_t({\bf n})$. 
We have thus defined the map $S:\Lambda_t\to\Gamma_t$, which we call
step map. Clearly, $S$ is invertible. The RSK dynamics on $\Lambda_t$
induces a dynamics of steps on $\Gamma_t$. By construction the RSK
dynamics stays inside $\Gamma_t$ with probability one.

Remarkably, the distribution of the multi-layer PNG model at time $t$
has a simple structure in being the uniform Lebesgue measure on all
admissible step configurations.

\begin{theorem}\label{the1}
Let $w_t$ be the uniform measure on $\Gamma_t$, which means that
$w_t \big(\Gamma_t(0)\big) = 1$ and $w_t \upharpoonright \Gamma_t
(\text{\bf n})$ is the $2|{\bf n}|$-dimensional Lebesgue measure on
$\Gamma_t (\text{\bf n})$. Then $\int w_t = Z(t) = \exp(2t^2)$ and
$\mu_t = Z(t)^{-1} w_t$ is a probability measure on
$\Gamma_t$. If the height lines evolve under the RSK dynamics, then
$\mu_t$ is the joint distribution of $\{h_\ell (x,t),\, x\in\IR,\,
\ell\in\IN_- \}$ under the step map $S$.   
\end{theorem}
{\it Proof}: 
The nucleation events determining the line configuration at time $t$
are a Poisson process of intensity $2$ in the 
triangle $\{(x',t')\big|\,\,|x'|<t',\,0\leq t'\leq t\}$. If
$(x',t')$ is a generic Poisson point, we label it through the new
coordinates $(y^+,y^-)=\big(x'-(t-t'),x'+(t-t')\big)$. A Poisson point
configuration consisting of $N$ points is then given by
$(y^+_i,y^-_i)_{i=1,\dots,N}$, such that $-t<y^+_1<\cdots<y^+_N<t$ and
$y^+_i<y^-_i$, $i=1,\dots,N$. The set of 
all such point configurations is denoted by $\Delta_t(N)$, considered
as a subset of $[-t,t]^{2N}$. We also set
$\Delta_t=\bigcup_{N\geq0}\Delta_t(N)$, with $\Delta_t(0)$ a single
point. $\Delta_t$ inherits from the Poisson process the probability
measure $\nu_t$, where $\nu_t\big(\Delta_t(0)\big)=e^{-2t^2}$ and
$\nu_t \upharpoonright \Delta_t(N)=e^{-2t^2}dy^+_1\cdots
dy^+_Ndy^-_1\cdots dy^-_N$.

Next we define the growth map $G:\Delta_t\to\Lambda_t$. For given
$(y^+_i,y^-_i)_{1\leq i\leq N}\in\Delta_t(N)$ we run the RSK dynamics
to determine the line configuration at time $t$. Conversely for a
given line configuration in $\Lambda_t$ we run the RSK dynamics
backwards in time which then determines the Poisson points
corresponding to the nucleation events. Thus $G$ is well defined. 
If all the line configurations with coinciding up- resp. down-step
positions in different lines---a set of Lebesgue measure zero under
the step map---are removed from  $\Lambda_t$, then $G$ is even bijective.

From the construction it is obvious that for $S\circ G:
(y^+_i,y^-_i)_{i=1,\dots,N}\mapsto\big((y^{\ell,+}_j,y^{\ell,-}_j)_{1\leq
  j\leq n_\ell}\big)_{\ell\leq0}$ with $N=|{\bf n}|$ one has the
following set equalities,
$\{y^+_i,\,i=1,\dots,N\}=\{y^{\ell,+}_j,\,j=1,\dots,n_\ell,\,\ell\in\IN_-\}$ 
and 
$\{y^-_i,\,i=1,\dots,N\}=\{y^{\ell,-}_j,\,j=1,\dots,n_\ell,\,\ell\in\IN_-\}$. 
Thus the map $S\circ G$ induces a mere relabeling of points. In
particular $\nu_t$ is transformed to $\mu_t$ under  $S\circ G$. 
$\Box$\medskip 

Equipped with Theorem \ref{the1} the reader may jump ahead to Section
\ref{sec.c} where the statistical properties of the line ensemble
$\mu_t$ are studied. We take a little detour to report on two
observations of interest. First there is a variant of the multi-layer
PNG dynamics which was introduced by Gates and Westcott \cite{Gat} in
modeling crystal growth, hence called GW dynamics by us. Bulk
properties of the GW dynamics are studied in 
\cite{Pra3} in the context of the anisotropic KPZ equation. 
Gates and Westcott regard the lines $h_\ell$ as contour lines of a
crystal surface. The crystal is made up of atomic two-dimensional
layers stacked along the $z$-axis. Layer $0$ and below are completely
filled. For layer $1$ only the domain $\{(x,y)\in\IR^2\big|\,y\leq
h_0(x,t)\}$ is filled with atoms, and in general, layer $-\ell+1$ is
filled in the domain $\{(x,y)\in\IR^2\big|\,y\leq h_\ell(x,t)-\ell\}$,
$\ell\in\IN_-$. 
The crystal is in contact with its supersaturated vapor. If overhangs
are not permitted, it is a natural modeling assumption that each
contour line grows under the PNG rules subject to a constraint of no
touching. In this sense the GW dynamics is ``more stochastic'' than
the RSK dynamics. Despite different rules, the RSK and GW dynamics
yield an identical distribution for the line configurations at time
$t$. Thus, alternatively, the proof of Theorem 
\ref{the1} could be based on the GW dynamics, cf.~below. 

A further observation is that the GW dynamics admits a discrete
space-time version which, as to be discussed at the end of this
section, inherits the simplicity of the distribution at time $t$. In
fact, the discrete version of the GW dynamics is isomorphic to the
shuffling algorithm for the Aztec diamond. This provides us with yet
another approach to Theorem \ref{the1}, namely to take the continuum
limit of its discrete analogue. 

Let us start with the GW dynamics. The top line $h_0$ evolves as the
PNG droplet. With nucleation rate $2$ pairs of up- and down-steps are
generated in the forward light cone of the origin and move apart with velocity
$\pm1$. Upon collision step pairs 
annihilate each other. The lower lines $h_{-1}$, $h_{-2}$, $\dots$
follow the same dynamics, independently of each other, under the
condition that nucleations are suppressed whenever they violate the 
monotonicity constraint, i.e.~for height line $\ell$, $\ell<0$,
nucleations occur only in the region
$\{x\in\IR\big|\,\,|x|<t,\,h_{\ell+1}(x,t)-h_\ell(x,t)\geq2\}$ with
space-time rate $2$. As for the RSK dynamics, we have 
to convince ourselves that the GW dynamics lives on $\Gamma_t$ with
probability one. The simple procedure is to use duality. We regard the
multi-layer PNG model as the stochastic evolution of the random field 
$\eta_t (j, x)$ over $\IZ \times \IR$ with values in $\{0,1\}$ by setting 
\begin{equation}
  \label{eq:new}
   \eta_t (j,x) = \left \{ 
    \begin{array}{ll} 
      1,\quad & \mbox{if $h_\ell (x,t) = j$ for some $\ell$},\\ 
      0, & \mathrm{otherwise}. \end{array}\right.
\end{equation}
The starting configuration is $\eta_0 (j,x) = 1$ for $j
\leq 0$, $\eta_0 (j,x)=0$ for $j \geq 1$. Let us introduce the flipped
configuration through $\bar{\eta}_t (j,x) = 1-\eta_t (j,x)$ with
the corresponding height lines $\bar{h}_\ell (x,t), \ell = 1,2, \dots$ .
Then $\bar{h}_1 (x,t)$ evolves again as the PNG droplet, now
growing downwards towards negative $j$. 
Thus for the original
process, at any given time $t$, there is a random index $\ell$ such
that $h_\ell (x,t) = \ell$ for all $x$. This proves that at any time $t$
with probability $1$ there are only finitely many up- and down-steps.
\begin{proposition}\label{prop1}
Let $\{h_\ell(x,t),\,x\in\IR,\,\ell\in\IN_-\}$ be the height lines as
generated by the GW dynamics. Under the step map $S$ their joint
distribution is $\mu_t$, $\mu_t$ of Theorem \ref{the1}.
\end{proposition}
{\it Proof}: We did not discover a global version comparable to the
proof of Theorem
\ref{the1} and have to rely on an infinitesimal argument. Let $L(t)$
be the forward generator of the Markov jump process induced on
$\Gamma_t$ through the GW dynamics. We have to prove that
\begin{equation}\label{b.a}
\frac{d}{dt} (e^{-2t^2}w_t) = L(t) e^{-2t^2} w_t
\end{equation}
which means
\begin{equation}\label{b.b}
- 4 t w_t = L(t) w_t \,,
\end{equation}
since $d w_t /dt = 0$. $L(t)$ has four pieces. (i): There is a loss
term $-4t$ due to nucleation events at $h_0(x,t)$ for $x\in[-t,t]$,
which cancels 
the left hand side of (\ref{b.b}). (ii): The free flow terms
$\partial/ \partial y^{\ell,+}_j$ of the up-step motion and
$-\partial/\partial y^{\ell,-}_j$ of the down-step motion vanish
when acting on the Lebesgue measure. There is no boundary
contribution at $\pm t$, because the boundary moves
with the same speed as the steps. 
(iii) and (iv): Let $r_\ell = | \{x\in\IR\big|\,\,|h_\ell (x,t)-
h_{\ell-1}(x,t)|\geq 2 \}|$ with 
$|\cdot|$ denoting the one--dimensional Lebesgue measure and let $r
= \sum_{\ell\leq0} r_\ell$. For given $|\text{\bf n}|$ the current
configuration gains in probability due to a transition from
$|\text{\bf n}|+1$ 
to $|\text{\bf n}|$ through the collision of an up-step and a
down-step. Their 
relative velocity is $2$. Thus for a small time interval $dt$ the gain
is $2r\,dt$, 
since the ratio of the weight at $|\text{\bf n}|+1$ to the weight at
$|\text{\bf n}|$ equals $1$. There is a loss of probability due to
nucleations in the 
current configuration. For the time $dt$ it is $2r\,dt$, since
nucleation events have intensity $2$. Thus the gain term (iii) cancels
exactly the loss term (iv). More extended versions of our argument can
be found in \cite{Gat,Pra3} where, however, a different geometry is
discussed.
$\Box$\medskip 

Next we introduce a discrete time version of the multi-layer 
PNG. As explained in \cite{Joh, Helf} this model is isomorphic to
the shuffling algorithm for the Aztec diamond. We discretize time,
now denoted by $\tau = 0,1, \dots$ . We also discretize the space axis
in units of $\delta$. As before, at time $\tau$, the height
lines are $h_\ell(x,\tau)$, $x \in \IR$, $\ell\in \IN$. $h_\ell(x,
\tau)=\ell$ for $|x| \geq \delta \tau$. The non-crossing constraint
$h_{\ell-1}(x,\tau) < h_\ell(x,\tau)$ is in force. Up- and down-steps
are allowed only at midpoints of the form $(m + \frac{1}{2})
\delta$, $m \in \IZ$ and their distances must be odd, i.e. of the
form $(2 m+1)\delta$. To update from time $\tau$ to time $\tau +1$
only changes inside the strip $[-(\tau +1)\delta, (\tau
+1)\delta]$ are allowed. The actual update consists of a
deterministic and a stochastic step.
\\
(i) \emph{Deterministic step:} Given $h_\ell(x,\tau)$ every up-step is
moved a $\delta$-unit to the 
left, every down-step a $\delta$-unit to the right. If at time
$\tau$ there is a block of length $2 \delta$, short a $2
\delta$-block, with a down-step to the left of an up-step, then they
annihilate each other, i.e. in this block $h_\ell(x, \tau)$ is
replaced by its maximum. The configuration after the
deterministic step is denoted by $\tilde{h}_\ell (x, \tau +1)$.
\\
(ii) \emph{Stochastic step:} The constant pieces of each height line
$\tilde{h}_\ell (x, \tau+1)$ are subdivided in consecutive $2
\delta$-blocks. To fix their location, the left endpoint of a $2
\delta$-block either coincides with the right endpoint of a $2
\delta$-block or is $\frac{1}{2} \delta$ away from an up-step, resp.
from a down-step.
If $\tilde{h}_0(x, \tau +1)= 0$, the $2 \delta$-blocks are of the
form $[(-\tau -1+2m)\delta, (-\tau -1+2m+2) \delta]$. If
$\tilde{h}_\ell (x, \tau +1)=\ell$ for all $x$ and $\tilde
h_{\ell+1}(x,\tau+1)\neq \ell+1$ for 
some $x$, the $2\delta$-blocks of $\tilde h_\ell(x,\tau+1)$ start at
$(y + \frac{1}{2}) \delta$ with 
$y$ the position of the first up-step (from the left) of
$\tilde{h}_{\ell+1}(x, \tau +1)$. Finally we disregard those $2
\delta$-blocks for which $\tilde{h}_{\ell+1}(x, \tau +1)- \tilde{h}_\ell
(x, \tau +1)=1$ for some $x$ inside the block. After these
preparations the stochastic update can be performed. Independently
for each $2 \delta$-block, we keep the original piece of the
height line with probability $1-q$, $0 < q < 1$, and otherwise nucleate an
up-step to the left and a down-step to the right midpoint of the two
adjacent $\delta$-intervals. The line configuration after the
stochastic update is denoted by $h_\ell(x,\tau+1)$. 

In the limit of rare events the discrete multi-layer PNG model
converges to its continuous time cousin. We set $t = \tau \delta$ and
denote by $[t]$ the integer part of $t$.
Then space-time is discretized in cells of lattice spacing
$\delta$. A nucleation event covers a block of two adjacent cells.
If we set $q = 4 \delta^2$, then in the limit $\delta \rightarrow
0$ we obtain a Poisson process of intensity $2$. Therefore in this
limit $h_\ell(x, [t/\delta]) \rightarrow h_\ell (x,t)$ as a stochastic
process.

The discretized multi-layer PNG model inherits the simplicity of
the time $\tau$ measure, $\mu(\tau)$. The height line
$h_\ell(x,\tau)$ has $n_\ell$ up-steps. The total number of up-steps
is then 
$\sum_{\ell\leq0} n_\ell = n$. To the collection of height lines
$\{h_\ell(x,\tau),\,x\in\IR,\,\ell\leq0\}$ we assign the weight
$(q/(1-q))^n$. The partition 
function, $Z_{\dd}(\tau)$, is the sum over all weights. We set
$Z_{\dd}(\tau)^{-1} = \IP\big(\{h_0(x,0)=0\}\big)$. Therefore
\begin{equation}\label{b.c}
Z_{\dd}(\tau) = \prod^{\tau}_{j=1}(1-q)^{-j} =
(1-q)^{- \tau(\tau +1)/2}\,.
\end{equation}
If the weight at time $\tau$ is denoted by $w(\tau)$,
we claim that
\begin{equation}\label{b.d}
\mu(\tau) = Z_{\dd}(\tau)^{-1} w(\tau)
\end{equation}
is the time  $\tau$ probability measure of the discrete multi-layer PNG.
Let $K_{\tau}$ be the transition kernel from $\tau$ to $\tau + 1$,
as explained in steps (i), (ii) above. We have to show $\mu(\tau +
1) = \mu(\tau) K_{\tau}$, equivalently
\begin{equation}\label{b.e}
(1-q)^{\tau +1} w(\tau +1) = w(\tau) K_{\tau}\,.
\end{equation}

(\ref{b.e}) is established in Proposition \ref{prop2} below. But
first we want to convince ourselves that $\mu(\tau)$ yields $\mu_t$ of
the continuous time PNG in the limit $\delta\to0$. We note that
for $\tau = [t/\delta], q = 4\delta^2$,
\begin{equation}\label{b.f}
\lim_{\delta \rightarrow 0} Z_{\dd}([t/\delta]) = \lim_{\delta
\rightarrow 0} \, \exp \big(- \textstyle\frac{1}{2}[t/\delta]([t/\delta]+1)
\log(1-4 \delta^2)\big) = e^{2t^2} = Z(t)\,.
\end{equation}
A configuration with $n$ up/down-step pairs has the weight
$(q/1-q)^n \cong(4 \delta^2)^n$. Because of the constraint in the
up-step locations, in the limit $\delta \rightarrow 0$ the weight
converges to the $2n$-dimensional Lebesgue measure constrained to
$\Gamma_t (\textbf{n})$, $n=|\text{\bf n}|$. Thus $\mu([t/\delta]) \to
\mu_t$ as $\delta \rightarrow 0$, as it should be.
\begin{proposition}\label{prop2}
Let the weight, $w(\tau)$, of the height lines of the discrete
multi-layer PNG be given by
$(q/(1-q))^n$, where $n$ is the total number of up-steps (equivalently
down-steps). Then (\ref{b.e}) holds.
\end{proposition}
{\it Proof}: Let $w(\tau+1)$ be the weight for the configuration
$h_\ell(x, \tau +1)$. We construct from it the configuration
$\tilde{h}_\ell (x, \tau +1)$ by removing all spikes from $h_\ell(x,
\tau +1)$, i.e. all $2\delta$-blocks containing an up-step to the
left and a down-step to the right. Let $s_\ell$ be the number of
spikes for $h_{\ell+1}(x, \tau+1)$, let $n_\ell$ be the number of
up-steps for $\tilde{h}_\ell(x, \tau+1)$, and let $b_\ell$ be the
number of $2\delta$-blocks with $\tilde h_\ell(x,\tau+1)$ constant such that
$\tilde{h}_{\ell+1}(x,\tau+1)-\tilde h_\ell(x,\tau+1)\geq 2$ within that
block, $b_0$ is the number of flat $2\delta$-blocks of $\tilde
h_0(x,\tau+1)$. Next we map the configuration $\tilde h_\ell(x,\tau+1)$
to the configuration $\tilde{h}_\ell(x, \tau)$ by moving all
up-down-steps one 
step backwards in time. By construction, $n_\ell$ does not change. Let
$a_\ell$ be the number of downwards open blocks, i.e.~flat $2\delta$-blocks
of $\tilde{h}_\ell(x, \tau)$ 
such that $h_{\ell-1} (x,\tau)$ has distance $\geq 2$ within that block.

The transition kernel $K_{\tau} = K_{\dd} K_{\ss}$, where $K_{\dd}$ is the
deterministic step (i) and $K_{\ss}$ is the stochastic step (ii). To
compute $w(\tau) K_{\tau}$, we first evaluate $w(\tau) K_{\dd}$ in the
configuration $\tilde{h}_\ell(x, \tau +1)$. We have to sum over all
line configurations $h_\ell(x,\tau)$ leading to $\tilde{h}_\ell(x,
\tau +1)$ under $K_{\dd}$. A downwards open block of
$\tilde h_\ell(x,\tau)$
had either no steps, weight 1, or a downwards spike, weight
$q/(1-q)$. Summing over these $2^{a_\ell}$ possibilities results in  
the weight $(1-q)^{-a_\ell}\big(q/(1-q)\big)^{n_\ell}$ for
$\tilde{h}_\ell(x,\tau+1)$. Applying  
the stochastic transition $K_{\ss}$ yields the weight $w(\tau)K_\tau$
evaluated at $h_\ell(x, \tau +1)$ as
\begin{equation}\label{b.g}
q^{s_\ell} (1-q)^{b_\ell - s_\ell}(1-q)^{-a_\ell} (q/(1-q))^{n_\ell}\,.
\end{equation}
On the other hand, according to $w(\tau +1)$, $h_\ell(x, \tau+1)$ has
the weight $\big(q/(1-q)\big)^{s_\ell+n_\ell}$. Comparing with (\ref{b.e}) and
(\ref{b.g}), we have to prove
\begin{equation}\label{b.h}
\sum_{\ell\leq0}(b_\ell - a_\ell) = \tau +1 \,.
\end{equation}

Let $N$ be the index of the last height line for which $h_N(x,
\tau)= N$ for all $x$. Then $a_N = 0$. For two adjacent lines it is
easily verified that
\begin{equation}\label{b.i}
b_\ell - a_{\ell+1} = n_{\ell+1}-n_\ell,\quad \ell\leq-1 \,.
\end{equation}
Inserting in the left side of (\ref{b.h}) and using $n_N = 0$, the
telescoping sum gives $b_0 + n_0$, which by definition is independent
of $\tilde{h}_0 (x, \tau +1)$ and equals $\tau + 1$. $\Box$

\section{$1+1$-dimensional Fermi field} \label{sec.c}
\setcounter{equation}{0}

$w_t$ is the uniform distribution on all allowed line
configurations of the continuous time multi-layer PNG. Except for
exclusion (entropic repulsion), the 
height lines do not interact. 
Such a statistical mechanics system is most conveniently analyzed
through the transfer matrix method. Its implementation requires the
height lines to be restricted to a bounded interval
$\{-M,-M+1,\dots,M\}=I_M$. The case of interest is then obtained in
the limit as $M\to\infty$. 
To explain the principle, we omit the
argument $t$ and label the height lines more conventionally as
$h_\ell(x)$, $\ell=1,\dots,N$, $N\leq 2M+1$, with $x\in[0,t]$.
The height lines
are constrained through $-M\leq h_1(x)<\dots<h_N(x)\leq M$ for all
$x\in[0,t]$. 
In addition we fix the initial configuration $q$ and the
final configuration $q'$, i.e. $h_\ell(0)=q_\ell$, $h_\ell(t)=q'_\ell$,
$\ell=1,\dots,N$. As before under the step map the $\ell$-th height line is
specified by the location of the up-steps,
$0<y_1^{\ell,+}<\cdots<y_{n_\ell}^{\ell,+}<t$, and the location of the
down-steps, $0<y_1^{\ell,-}<\cdots<y_{n'_\ell}^{\ell,-}<t$, where
$n_\ell\neq n'_\ell$ is allowed. Admissible line configurations
are assumed to have a 
uniform weight, which means that a small volume element has the weight
$\prod_{\ell=1}^N
\prod_{j=1}^{n_\ell}\prod_{j'=1}^{n'_\ell}dy^{\ell,+}_jdy^{\ell,-}_{j'}$.
The configuration with no steps has weight $1$. 

We want to compute the partition function $Z_t(q,q')$ which is
defined as the weight integrated over all admissible step
configurations. For this purpose the simplex
$\Omega_N=\{q\in\IZ^N\big|\,-M\leq q_1<\cdots<q_N\leq M\}$ is
introduced. Clearly 
$q,q'\in\Omega_N$ and we regard $Z_t(q,q')$ as a
$|\Omega_N|\times|\Omega_N|$ matrix. By the product property of the
Lebesgue measure $Z_t$ satisfies the semigroup property
\begin{equation}
  \label{eq:semigroup}
  Z_tZ_s=Z_{t+s},\quad t,s\geq 0,\,\,Z_0=\id,
\end{equation}
$\id$ the identity matrix.
Thus there exists an infinitesimal generator $G_N$, such that
\begin{equation}
  \label{eq:infgen}
  Z_t=e^{-tG_N},\quad t\geq0.
\end{equation}
Differentiating at $t=0$ one concludes that $G_N$ acting on functions
$f$ on $\Omega_N$ is given by
\begin{equation}
  \label{eq:generator}
  G_Nf(q)=-\sum_{q'\in\Omega_N}c(q,q')f(q'),
\end{equation}
where for $q,q'\in\Omega_N$
\begin{equation}
  \label{eq:rates}
  c(q,q')=\left\{
    \begin{array}{cl}
      1&\mbox{if $\sum_{\ell=1}^N|q_\ell-q'_\ell|=1$,}\\
      0&\mbox{otherwise}.
    \end{array}
  \right.
\end{equation}

Computationally much more powerful is to impose the constraint of no
overlap through antisymmetry. Let ${\cal F}_N$ be the subspace of
$\ell_2\big((I_M)^N\big)$ consisting of antisymmetric functions over
$(I_M)^N$, i.e. $f\in{\cal F}_N$ satisfies
\begin{equation}
  \label{eq:antisym}
  f(q_1,\dots,q_N)=(-1)^{\text{sign}\,\pi}f(q_{\pi(1)},\dots,q_{\pi(N)})
\end{equation}
for every permutation $\pi$. ${\cal F}_N$ is equipped with the
canonical basis $f_q, q\in\Omega_N$, defined through
\begin{equation}
  \label{eq:basis}
  f_q(q')=\frac1{\sqrt{N!}}\sum_\pi(-1)^{\text{sign}\,\pi}
  \delta_q(q'_{\pi(1)},\dots,q'_{\pi(N)}) 
\end{equation}
with $\delta_q(q')=1$ if $q=q'$ and $\delta_q(q')=0$ otherwise. The
normalization is chosen such that $\langle
f_q,f_{q'}\rangle=\delta_q(q')$, with $\langle\cdot,\cdot\rangle$
denoting the scalar product in ${\cal F}_N$.
Let us also define the one-particle
Hamiltonian
\begin{eqnarray}
  \label{eq:onepartham}
  H^M_{\text{d}}\psi(-M)=-\psi(-M+1),\quad
  H^M_{\text{d}}\psi(M)=-\psi(M-1),
  \nonumber\\
  H^M_{\text{d}}\psi(n)=-\psi(n+1)-\psi(n-1)\quad\mbox{for $|n|<M$},
\end{eqnarray}
as  acting on functions $\psi$ over $I_M$. The corresponding $N$-particle
Hamiltonian, is then given through
\begin{equation}
  \label{eq:npartham}
  H^M_{\text{d},N}=\sum_{j=1}^N1\otimes\cdots\otimes H^M_{\text{d}}
  \otimes\cdots\otimes1,
\end{equation}
where $H^M_{\text{d}}$ is inserted at the $j$-th position of the
$N$-fold product. Clearly for  $f\in{\cal F}_N$ one has
$H^M_{\text{d},N}f\in{\cal F}_N$ and $H^M_{\text{d},N}$ is regarded
as acting on ${\cal F}_N$. With these notations one has the identity
\begin{equation}
  \label{eq:fockeq}
  \langle f_q,e^{-tH^M_{\text{d},N}}f_{q'}\rangle=e^{-tG_N}(q,q')
  =Z_t(q,q')
\end{equation}
for $q,q'\in\Omega_N$.

At this point it is more convenient to switch to fermionic language
which is devised precisely to take the antisymmetry into account. The
CAR algebra over $I_M$ is generated by $a^*(j)$, $a(j)$, $j\in I_M$.
They satisfy the canonical anticommutation relations 
\begin{equation}
  \label{eq:CARrel}
  \{a(i),a^*(j)\}=\delta_{ij},\quad
  \{a(i),a(j)\}=0,\quad\{a^*(i),a^*(j)\}=0,
\end{equation}
$i,j=-M,\dots,M$, $\{A,B\}=AB+BA$. In the Fock representation the
algebra is realized as operators on the antisymmetric Fock space
${\cal F}$ over $I_M$,
\begin{equation}
  \label{eq:Fockspace}
  {\cal F}=\bigoplus_{N=0}^{2M+1}{\cal F}_N.
\end{equation}
The second quantization of the $(2M+1)\times(2M+1)$ matrix
$H^M_{\text{d}}$ is defined by
\begin{equation}
  \label{eq:secondquant}
  \widehat{H}^M_{\text{d}}=\sum_{i,j\in I_M}a^*(i)
  \big(H^M_{\text{d}}\big)_{ij}a(j).
\end{equation}
$\widehat{H}^M_{\text{d}}$ restricted to ${\cal F}_N$ agrees with
$H^M_{\text{d},N}$. From (\ref{eq:fockeq})
one concludes that as a fermionic operator the transfer matrix is
given through  
\begin{equation}
  \label{eq:transfermatrix}
  e^{-t\widehat{H}^M_{\text{d}}},\quad t\geq0,
\end{equation}
which covers all $0\leq N\leq 2M+1$.

We exploit the new flexibility by assuming that
$q,q'\in\Omega=\bigcup_{N=0}^{2M+1}\Omega_N$, which is identified with
$\{0,1\}^{I_M}$. The case of interest is $q=q'$ and each boundary
configuration $q$ has the weight
$\prod_{j=-M}^M\exp\big(\lambda(j)\sum_{i=1}^N\delta_{q_i,j}\big)$.
In other words the configuration $q$ has a product weight with factor
$e^{\lambda(j)}$ if site $j$ is occupied and factor $1$ if the site
$j$ is empty.
The corresponding partition function is then given through
\begin{equation}
  \label{eq:partfunc}
  Z^\lambda_M=\tr\big[ e^{\widehat{N}^M}
  e^{-t\widehat{H}^M_{\text{d}}}\big],
\end{equation}
where the trace is over ${\cal F}$ and 
\begin{equation}
  \label{eq:numberop}
  \widehat{N}^M=\sum_{j\in I_M}\lambda(j)a^*(j)a(j).
\end{equation}
The probability that site $j$ at $x$ is occupied, i.e. $h_\ell(x)=j$
for some $\ell$, is obtained from the transfer matrix as
\begin{equation}
  \label{eq:occupprob}
  \big(Z^\lambda_M\big)^{-1}\tr\big[e^{\widehat{N}^M}
  e^{-x\widehat{H}^M_{\text{d}}}a^*(j)a(j)
  e^{-(t-x)\widehat{H}^M_{\text{d}}}\big]=\tr[\widehat{\rho}_x^Ma^*(j)a(j)],
\end{equation}
with the density matrix
\begin{equation}
  \label{eq:densityop}
  \widehat{\rho}_x^M=\big(Z^\lambda_M\big)^{-1}
  e^{-(t-x)\widehat{H}^M_{\text{d}}}e^{\widehat{N}^M}
  e^{-x\widehat{H}^M_{\text{d}}},
\end{equation}
$0\leq x\leq t$.

Exponentials of operators quadratic in $a,a^*$ are easily handled. If
$A$ is a $(2M+1)\times(2M+1)$ matrix with second quantization
\begin{equation}
  \label{eq:secondqu}
  \widehat{A}=\sum_{i,j\in I_M}a^*(i)A_{ij}a(j),
\end{equation}
then 
\begin{equation}
  \label{eq:nparticle}
  e^{\widehat{A}}=e^A\otimes\cdots\otimes e^A
\end{equation}
on ${\cal F}_N$. This implies
\begin{equation}
  \label{eq:trexp}
  \tr[e^{\widehat{A}}]=\det(1+e^A)=Z_A,
\end{equation}
compare with \cite{Reed}. Let us set
$\widehat{\rho}_A=(Z_A)^{-1}e^{\widehat{A}}$ as density matrix. The
two-point function has the form 
\begin{equation}
  \label{eq:twopointfunction}
  R(i,j)=\tr[\widehat{\rho}_A a^*(i)a(j)]=\Big((1+e^{-A})^{-1}\Big)_{ij}
\end{equation}
and more generally
\begin{equation}
  \label{eq:npointfunction}
  \tr[\widehat{\rho}_A a^*(i_1)\cdots a^*(i_m)a(j_n)\cdots a(j_1)]=
  \delta_{m,n}\det\big(R(i_k,j_{k'})\big)_{k,k'=1,\dots,m}.
\end{equation}
The expectations of other monomials are determined by means of the
anticommutation relations (\ref{eq:CARrel}). 
One may regard
$\tr[\widehat{\rho}_A\,\cdot\,]=\omega_A(\,\cdot\,)$ as a linear
functional on the CAR algebra. By definition $\omega_A(1)=1$. If in
addition $\omega_A$ is positive, then $\omega_A$
is called a quasifree state \cite{Brat}.

As can be seen from (\ref{eq:nparticle}) products of exponentials
follow the same pattern. For the $(2M+1)\times(2M+1)$ matrices $A,B$
we set 
\begin{equation}
  \label{eq:logc}
  e^Ae^B=e^C.
\end{equation}
Then
\begin{equation}
  \label{eq:logc2nd}
  e^{\widehat{A}}e^{\widehat{B}}=e^{\widehat{C}}
\end{equation}
with $\widehat{\,\cdot\,}$ defined as in (\ref{eq:secondqu}).

Each $a^*(i)a(i)$ is a symmetric projection and thus has eigenvalues
in $\{0,1\}$. $\{a^*(i)a(i),\,i\in 
I_M\}$ is a family of commuting operators. Under $\widehat{\rho}^M_x$
they have a joint signed spectral measure which by construction is a
probability measure on $\{0,1\}^{I_M}$. By
(\ref{eq:npointfunction}) it is of determinantal form, compare also with 
(\ref{c.r}) below. Thus under $\widehat{\rho}^M_x$ the family
$\{a^*(i)a(i),\,i\in I_M\}$ is a point process on $I_M$. In the
probabilistic literature 
point processes of this structure are known as determinantal \cite{Sosh}.

Our construction extends to occupation variables depending on several
$x$. For example, the probability that site $j$ at $x$ and site $i$
at $y$, $0<x<y<t$, are both occupied is obtained through the transfer
matrix as
\begin{equation}
  \label{eq:twoptcorr}
  \big(Z^\lambda_M\big)^{-1}\tr\big[
  e^{\widehat{N}^M}e^{-x\widehat{H}^M_{\text{d}}}a^*(j)a(j)
  e^{-(y-x)\widehat{H}^M_{\text{d}}}a^*(i)a(i)e^{-(t-y)
  \widehat{H}^M_{\text{d}}}\big].    
\end{equation}
Such expressions can be computed using (\ref{eq:logc2nd}) and
(\ref{eq:logc}), and the rules
\begin{equation}
  \label{eq:anotherformula}
  e^{t\widehat{A}}a(j)e^{-t\widehat{A}}=\sum_{i\in I_M}(e^{-tA})_{ji}a(i),
  \quad
  e^{t\widehat{A}}a^*(j)e^{-t\widehat{A}}=\sum_{i\in I_M}a^*(i)(e^{tA})_{ij}.
\end{equation}

With these preparations we return to the PNG droplet. Recall the
definition (\ref{eq:new}) of the random field of  
occupation variables $\eta_t(j,x)$, $j \in \IZ$, $x \in \IR$: $\eta_t 
(j,x) = 1$ if $h_\ell (x,t) =j$ for some $\ell$ and $\eta_t (j,x) =0$
otherwise. By definition $\eta_t(j,x) = 0$ for $j \geq 1$, $|x|
\geq t$, and $\eta_t(j,x) =1$ for $j \leq 0$, $|x| \geq t$. The joint
distribution of $\eta_t(j,x)$, $j\in\IZ$, $x\in\IR$, is induced through the
probability measure $\mu_t$ with expectation denoted by $\IE_t$. 
Our first goal is to obtain the joint distribution of $\eta_t(j,0)$,
$j\in\IZ$. From the considerations above, it is obvious that this
point process is determinantal. The only remaining task is to compute
the two-point function and to study its limit behavior. One should pay
attention to a minor linguistic problem. The meaning of the time
parameter $t$ of the PNG model is now reduced to a mere scaling
parameter. $j$ labels fermionic space and $x$ stands for fermionic
time. Thus space-time means now $\IZ\times\IR$.

Let us start by explaining the result for moments at $x=0$.
We introduce the limit $M\to\infty$ of the one-particle Hamiltonian 
in (\ref{eq:onepartham}) as
\begin{equation}
  \label{eq:Hd}
  H_{\text{d}}\psi(n)=-\psi(n+1)-\psi(n-1).
\end{equation}
In addition, we add a linear potential of slope $1/t$ to
define 
\begin{equation}\label{c.c}
H_t \psi (n) = - \psi(n+1) - \psi (n-1) + \frac{n}{t} \psi (n),
\end{equation}
regarded as an operator on $\ell_2=\ell_2(\IZ)$. $H_t$ has the
complete set of 
eigenfunctions $\varphi^{(l)}(n)=J_{n-l}(2t)$, $l\in\IZ$, with
eigenvalues $\varepsilon_l=\frac{l}{t}$,
$H_t\varphi^{(l)}=\varepsilon_l\varphi^{(l)}$,
$\langle\varphi^{(l)},\varphi^{(l')}\rangle=\delta_{l\,l'}$,
$\langle\cdot,\cdot\rangle$ denoting the scalar product. Here
$J_n(z)$ is the Bessel function of integer order $n$ and we follow throughout
the conventions of \cite{Abra}, Chapter 9.
We will need the spectral projection $B_t$ onto $\{H_t \leq 0\}$. In
position space its integral kernel is the discrete Bessel kernel
\begin{equation}\label{c.d}
B_t(i,j) = \sum_{l\leq0} J_{i-l}(2t) J_{j-l}(2t).
\end{equation}
Using that $H_t\, \varphi^{(l)} = \varepsilon_l\, \varphi^{(l)}$,
(\ref{c.d}) can be converted into a telescoping sum with the
result
\begin{equation}\label{c.e}
B_t (i,j) = \frac{t}{i-j} \big( J_{i-1}(2t) J_j(2t) -
J_i(2t)J_{j-1}(2t)\big)
\end{equation}
for $i \neq j$ and on the diagonal
\begin{equation}
  \label{eq:c.e1}
  B_t(i,i)=t \big( L_{i-1}(2t) J_i(2t) -
L_i(2t)J_{i-1}(2t)\big),
\end{equation}
where $L_j(2t)=\frac{d}{dj}J_j(2t)$.
We also introduce the CAR algebra ${\cal A}_{\dd}$ over $\IZ$. It
is generated by the operators $a(j),a^*(j),\,j\in\IZ$, satisfying the
canonical anticommutation relations (\ref{eq:CARrel}).  
Let $\omega_t$ be the quasifree state on the CAR algebra ${\cal
  A}_{\dd}$ defined through
$\omega_t\big(a(j)\big) = 0=\omega_t\big(a^*(j)\big)$ and the
two-point function
\begin{equation}\label{c.f}
\omega_t\big(a^*(i)a(j)\big) = B_t(i,j),
\end{equation}
which means that higher order monomials satisfy
(\ref{eq:npointfunction}) with $R$ replaced by $B_t$ \cite{Brat}. 
$\omega_t$ is the ground state for non-interacting fermions with
one-particle Hamiltonian (\ref{c.c}) at zero chemical potential. 
\begin{theorem}\label{the3}
We have
\begin{equation}\label{c.g}
\IE_t \bigg(\prod^m_{k=1} \eta_t(j_k, 0)\bigg) = \omega_t
\bigg(\prod^m_{k=1} a^*(j_k)a(j_k) \bigg)
= \det \big(B_t (j_k, j_{k'})\big)_{1\leq k,k'\leq m}\,,
\end{equation}
the second equality being valid only for pairwise distinct points
$j_1, \dots, j_m$.
\end{theorem}
{\it Proof}: The state $\omega_t$ has an infinite number of fermions
and cannot 
be represented as a vector in Fock space. Thus we first have to
constrain to finite volume $I_M$, such that $|h_\ell(x,t)| \leq
M, \ell=0, \dots, -M$. The corresponding uniform distribution is denoted
by $\mu^M_t$. Clearly  $\mu^M_t$ converges to $\mu_t$ as
$M\to\infty$. By construction $h_\ell(t,t)=h_\ell(-t,t)=\ell$,
$\ell=0,-1,-2,\dots\,$ . In other words $\eta_t(j,-t)=\eta_t(j,t)=1$ for
$j\leq0$ and $=0$ for $j\geq1$. It is convenient to approximate these
boundary configurations through the weight $e^\beta$, $-M\leq j\leq0$, weight
$e^{-\beta}$, $1\leq j\leq M$, for an occupied site and weight $1$ for an
empty site in the limit $\beta\to\infty$. Thus, if we set
$\widehat{N}^M$ as in (\ref{eq:numberop}) with
$\lambda(j)=1$ for $-M\leq j\leq0$ and $\lambda(j)=-1$ for $1\leq
j\leq M$, we have
\begin{equation}\label{c.h}
\IE^M_t \bigg(\prod^m_{k=1} \eta_t(j_k, 0) \bigg) = \lim_{\beta
\rightarrow  \infty} \frac{1}{Z(\beta)} \mathrm{tr} [e^{\beta
\widehat{N}^M} e^{-t \widehat{H}^M_{\dd}} \prod^m_{k=1} a^*(j_k)
a(j_k) e^{-t \widehat{H}^M_{\dd}} ]
\end{equation}
with $Z(\beta)$ the normalizing partition function. Since
the moments of a quasifree state are determined by the 
two-point function, to prove 
(\ref{c.g}) it suffices to consider the expectation of $a^*(j)
a(i)$ and to subsequently take the limit $M\rightarrow \infty$.
Let $P^M_-$ be the projection onto $\{-M,
\dots, 0 \}$, $P^M_+$ onto $\{1, \dots, M \}$, $P^M_+ +P^M_-$ being the
identity. Then 
\begin{eqnarray}\label{c.i}
\lefteqn{\lim_{\beta \rightarrow \infty} Z(\beta)^{-1}
  \mathrm{tr}
\big[e^{\beta \widehat{N}^M} e^{- t \widehat{H}^M_{\dd}}
a^*(j) a(i) e^{-t \widehat{H}^M_{\dd}}\big]}
\nonumber \\
&=& \lim_{\beta \rightarrow \infty} \Big(\big(1 + e^{t H^M_{\dd}}(e^{\beta}
P^M_+ + e^{-\beta} P^M_- )e^{t H^M_{\dd}}\big)^{-1}\Big)_{ij} 
\nonumber \\
&=& \Big(e^{-t H^M_{\dd}} P^M_- (P^M_++P^M_- e^{-2tH^M_{\dd}}
  P^M_-)^{-1} P^M_- e^{-tH^M_{\dd}} \Big)_{ij}\,.
\end{eqnarray}
If $P_-$ denotes the projection
onto $\IN_-$, $P_+=1-P_-$, then 
\begin{equation}
  \label{c.i1}
  \lim_{M\to\infty}P^M_\pm=P_\pm,\qquad
  \lim_{M\to\infty}e^{tH_{\text{d}}^M}=e^{tH_{\text{d}}}.
\end{equation}
To prove the theorem we only have to check the identity
\begin{equation}\label{c.j}
e^{-tH_{\dd}} P_- (P_++P_- e^{-2t H_{\dd}}P_-)^{-1} P_- e^{- t
  H_{\dd}} = B_t 
\end{equation}
as an operator identity on $\ell_2$.

We define the left shift $D$, $D \psi (n) = \psi(n+1)$, and the adjoint
right shift $D^*$, $D^* \psi (n) = \psi (n-1)$. Clearly $[D, D^*] = 0$.
One has $H_{\dd} = -D -D^*$. Using $2 \frac{d}{dt}J_n(t) = J_{n-1}(t)
- J_{n+1} (t)$, one obtains
\begin{equation}\label{c.k}
\frac{d}{dt} B_t = (D^* - D) B_t - B_t(D^* - D)\,.
\end{equation}
Integrating with the initial condition $B_0 = P_-$ yields
\begin{equation}\label{c.l}
B_t = e^{t(D^* - D)} P_- e^{-t (D^* - D)}
\end{equation}
and
\begin{equation}\label{c.m}
e^{tH_{\dd}}B_t e^{t H_{\dd}} = e^{-t(D+D^*)} e^{t(D^* - D)} P_- e^{-t(D^* -
D)} e^{-t (D+D^*)} = e^{-2tD} P_- e^{-2tD^*}\,.
\end{equation}
Therefore (\ref{c.j}) is equivalent to
\begin{equation}\label{c.n}
e^{-2tD} P_- e^{-2 tD^*} = (P_- e^{2 tD^*} e^{2 tD} P_-)^{-1}
\end{equation}
as an operator identity on $P_- \ell_2$. We decompose our space as
$\ell_2 = P_+ \ell_2 \oplus P_- \ell_2$. Then with the definition
\begin{equation}\label{c.o}
e^{-2tD} = {a \quad 0 \choose b \quad c}
\end{equation}
we have
\begin{equation}\label{c.p}
e^{-2tD} P_- e^{-2tD^*} = {a \quad 0 \choose b \quad c}{0 \quad 0
\choose 0 \quad 1} {a^* \quad b^* \choose 0 \quad c^*} = 
\bigg(
\begin{matrix}
  0&0\\0&c^* c\\
\end{matrix}
\bigg)
\,.
\end{equation}
Using the splitting of (\ref{c.o}), one constructs the inverse
operators $e^{2 t D}, e^{2 t D^*}$. By a straightforward
computation one obtains
\begin{equation}\label{c.q}
P_- e^{2 t D^*} e^{2 t D} P_- =
\bigg(
  \begin{matrix}
    0&0\\0&(c^* c)^{-1}
  \end{matrix}
\bigg).\quad\Box
\end{equation}

From Theorem \ref{the3} one immediately infers the distribution of the
height of the PNG droplet at $x=0$. Clearly
\begin{eqnarray}
  \label{eq:heightat0}
\IP_t\big(\{h(0,t)<n\}\big)&=&
  \IP_t\big(\{\eta_t(j,0)=0 \mbox{ for all $j\geq n$}\}\big)
  \nonumber\\
  &=&
\lim_{\beta\to\infty}\omega_t\bigg(\prod_{j=n}^\infty e^{-\beta
    a^*(j)a(j)}\bigg)
  \,\,=\,\,\lim_{\beta\to\infty}\det\big(1-(1-e^{-\beta})P_nB_t\big)
  \nonumber\\
  &=&
\det(1-P_nB_t),
\end{eqnarray}
where $P_n$ denotes the projection onto $\{n,n+1,\dots\}$ in
$\ell_2$. Since $L(0,t)=h(0,t)$, see eq.~(\ref{a.h}), we have
rederived 
that the length of the longest increasing subsequence of a Poissonized
random permutation has a distribution linked to the discrete Bessel
kernel. Previous proofs take the route via the Plancherel measure. We
refer to \cite{Johplanch}. It would be of interest to better
understand how these proofs are linked to the multi-layer PNG. 

So far we considered only the distribution of $\eta_t(j,x)$ at
$x=0$.
The transfer matrix method can handle also the distribution
referring to several $x$, like the joint distribution of
$\eta_t(i,0)$, $\eta_t(j,x)$, see eq.~(\ref{eq:twoptcorr}). The
transfer matrix is generated by the 
Hamiltonian $H_{\dd}$ of (\ref{eq:Hd}).
As in the case of fixed $x$, the joint moments have determinantal form
with 
the entries given by the Euclidean Fermi propagator $B_t(j,x;j',x')$. 
Following the scheme in (\ref{c.h}) it is defined through a finite volume approximation,
$B_t(j,x;j',x')=\lim_{M\to\infty}\lim_{\beta\to\infty}
B_t^{M\beta}(j,x;j',x')$, where 
\begin{equation}
  \label{eq:twoptcorrdyn1}
  B_t^{M\beta}(j,x;j',x')=\left\{
    \begin{array}{c}
      Z(\beta)^{-1}\tr\big[e^{\beta\widehat{N}^M}
      e^{-t\widehat{H}^M_{\text{d}}}
      \big(e^{-x\widehat{H}^M_{\text{d}}}a^*(j)
      e^{x\widehat{H}^M_{\text{d}}}\big)
      \big(e^{-x'\widehat{H}^M_{\text{d}}}a(j')
      e^{x'\widehat{H}^M_{\text{d}}}\big)
      e^{-t\widehat{H}^M_{\text{d}}}\big]
      \\[1mm]
      \hfill\mbox{for $-t\leq x\leq x'\leq
        t$,}
      \\[2mm]
      -Z(\beta)^{-1}\tr\big[e^{\beta\widehat{N}^M}
      e^{-t\widehat{H}^M_{\text{d}}}    
      \big(e^{-x'\widehat{H}^M_{\text{d}}}a(j')
      e^{x'\widehat{H}^M_{\text{d}}}\big)
      \big(e^{-x\widehat{H}^M_{\text{d}}}a^*(j)
      e^{x\widehat{H}^M_{\text{d}}}\big)
      e^{-t\widehat{H}^M_{\text{d}}}\big]
      \\[1mm]
      \hfill\mbox{for $-t\leq x'<x\leq
        t$}.    
    \end{array}\right.
\end{equation}
Note that time order must be respected in such a way that there are
only decaying exponentials. The minus sign in (\ref{eq:twoptcorrdyn1})
for $x'<x$
results from commuting $a^*$ and $a$. 
At coinciding arguments the definition conforms with
$\IE\big(\eta_t(j,x)\big)=B_t(j,x;j,x)$. 
Using (\ref{eq:anotherformula}) one obtains
\begin{equation}
  \label{eq:prop2}
   B_t(j,x;j',x')=\big(e^{-xH_{\dd}}(B_t-\id\theta(x-x'))
   e^{x'H_{\dd}}\big)_{jj'}, 
\end{equation}
for $|x|\leq t$, $|x'|\leq t$, $x\neq x'$,
with the step function $\theta(x)=0$ for $x<0$, $\theta (x)=1$ for
$x>0$. $B_t$ has a jump discontinuity at $x=x'$. 
For coinciding time arguments one has 
\begin{equation}
  \label{eq:prop2a}
  B_t(j,x;j',x)=(e^{-x H_{\dd}}B_te^{x H_{\dd}})_{jj'}.
\end{equation}

For later use the
propagator is rewritten in the eigenbasis of $H_t$. 
The integer order Bessel 
function has the representation 
\begin{equation}\label{c.aa}
J_n(2t)  = \frac{1}{2 \pi i} \oint\frac{dz}{z} e^{t(z^{-1}-z)}
z^n \,
\end{equation}
where the contour integration is a circle around $z = 0$.
Therefore
\begin{equation}\label{c.ab}
\big(e^{-x H_{\dd}} J_. (2t)\big)_n = \frac{1}{2 \pi i} \oint
\frac{dz}{z} e^{t(z^{-1}-z)} e^{x(z^{-1}+z)} z^n\,.
\end{equation}
Substituting $z$ by $(t+x)^{1/2}(t-x)^{1/2}z$ yields
\begin{equation}\label{c.ac}
\big(e^{-xH_{\dd}}J_.(2t)\big)_n = J_n \big(2 \sqrt{t^2 - x^2}\big)
\Big(\frac{t+x}{t-x}\Big)^{n/2}
\end{equation}
and thus for $x\neq x'$
\begin{eqnarray}
  \label{eq:prop3}
  B_t(j,x;j',x')&=&\sum_{l\in\IZ}\text{sgn}(x'-x)\theta\big((x-x')
  (l+{\textstyle\frac12})\big)\Big(\frac{t+x}{t-x}\Big)^{(j-l)/2}
  J_{j-l}\big(2\sqrt{t^2-x^2}\big)
  \nonumber\\
  &&\times J_{j'-l}\big(2\sqrt{t^2-{x'}^2}\big)
  \Big(\frac{t-x'}{t+x'}\Big)^{(j'-l)/2}.
\end{eqnarray}
At coinciding arguments  $x=x'$ one has 
\begin{equation}
  \label{eq:prop4}
  B_t(j,x;j',x)=\sum_{l\leq0}\Big(\frac{t+x}{t-x}\Big)^{(j-j')/2}
    J_{j-l}\big(2\sqrt{t^2-x^2}\big)J_{j'-l}\big(2\sqrt{t^2-{x}^2}\big)
\end{equation}
With these preparations for a general moment of the 
density field $\eta_t(j,x)$ one has the identity
\begin{eqnarray}\label{c.r}
  \lefteqn{\IE_t \bigg(\prod^m_{k=1} \eta_t(j_k, x_k)\bigg)=}
  \nonumber\\
  &=&\lim_{M\to\infty}\lim_{\beta\to\infty}Z(\beta)^{-1}
  \text{tr}\Big[e^{\beta\widehat{N}^M} 
  e^{-t\widehat{H}_{\dd}^M}\bigg(\prod_{k=1}^m e^{-x_{\pi(k)}
  \widehat{H}_{\dd}^M} a^*(j_{\pi(k)}) a(j_{\pi(k)}) 
  e^{x_{\pi(k)} \widehat{H}_{\dd}^M}\bigg) e^{-t \widehat{H}_{\dd}^M}\Big]
  \nonumber\\
  &=&\det\Big(B_t(j_k,x_k;j_{k'},x_{k'})\Big)_{1\leq k,k'\leq m}.
\end{eqnarray}
As written, (\ref{c.r}) is valid only for pairwise distinct
$x_1,\dots,x_m$, where $\pi$ is the unique 
permutation of $1,\dots,m$ such that the time ordering $-t\leq
x_{\pi(1)}<\cdots<x_{\pi(m)}\leq t$ is ensured. The spatial arguments
$j_1,\dots,j_m$ are arbitrary.
 While each off-diagonal factor in the determinant has a 
jump discontinuity at $x_k=x_{k'}$, the determinant itself depends
continuously on $x_1,\dots,x_m$, and thereby the continuous extension of
(\ref{c.r}) holds for all $j_k\in\IZ$, $-t\leq x_k\leq t$, $k=1,\dots,m$.

As an application of (\ref{c.r}) we establish the joint
distribution of $\{\eta_t(j,x), j \in \IZ \}$ for fixed $x$, $|x|
\leq t$. From (\ref{eq:prop4}) one derives immediately
\begin{equation}\label{c.ad}
e^{-x H_{\dd}} B_t e^{x H_{\dd}} = g \,B_{\sqrt{t^2-x^2}}\,g^{-1}
\end{equation}
where $g$ is a multiplication operator with diagonal
entries $g(n) =\big((t+x)/(t-x)\big)^{n/2}$.
In (\ref{c.r}) we take the limit of
coinciding $x_k=x$, $k=1,\dots,m$, leaving $j_1,\dots,j_m$ pairwise
distinct. Upon forming the determinant in
(\ref{c.r}) the similarity transformation $g$ drops out and the result
is (\ref{c.g}) with
$B_t$ replaced by $B_{\sqrt{t^2-x^2}}$. Thus the joint distribution of
$\{\eta_t(j,x), j \in \IZ \}$ is again given through the discrete
Bessel kernel with time parameter modified from $t$ to
$\sqrt{t^2-x^2}$.   

The same conclusion
can be drawn by taking the discrete PNG model as starting point. As
explained in \cite{Joh} the analogue of the fixed $x$
distributions is given through the Krawtchouk polynomials. Their
limit as $\delta \rightarrow 0$, $q = 4 \delta^2$, $\delta \tau =
t$, yields the joint distribution of $\{\eta_t(j,x), j \in \IZ \}$
as given through (\ref{c.g}) with parameter $\sqrt{t^2 - x^2}$ instead
of $t$.

In the following section we establish the scaling limit of the PNG
droplet at locations of order
$(yt^{2/3},2t+ut^{1/3})$, $y,u\in\IR$.
Since at $x=wt$, $|w|< 1$, the distribution is
determined by the 
discrete Bessel kernel $B_{t \sqrt{1-w^2}}$, we could instead of $w = 0$
choose any other reference point $(wt,2\sqrt{1-w^2}t)$ and relative
displacements 
$(wt + yt^{2/3},2\sqrt{1-w^2}t+ut^{1/3})$. Except for scale factors,
the limit $t 
\rightarrow \infty$ does not depend on the choice of $w$.

\section{Edge scaling, convergence to the Airy process} \label{sec.d}
\setcounter{equation}{0}
We plan to establish that the statistics of the PNG droplet close to
$x=0$ converges to the Airy process. Since only moments are under
control, the natural strategy is to prove that the Fermi field of
Section \ref{sec.c} has a limit when viewed from the density edge. In
particular, this implies, that the statistics of the last fermionic
world line has a limit, which is the desired result.

Since $\langle h(0,t)\rangle=2t$, the focus has to be at
$x=0+yt^\alpha$, $j=2t+u t^\beta$, 
$\,\,y,u\in\IR$ fixed. Rescaling
$\partial\psi/\partial x=H_t\psi$ accordingly one obtains
\begin{equation}
  \label{eq:1}
  \frac{\partial}{\partial y}\psi=
  t^\alpha\big(-\psi(u+t^{-\beta})-\psi(u-t^{-\beta})+
  \frac1t(2t+ut^\beta)\psi(u)\big).
\end{equation}
To have a limit operator as $t\to\infty$ requires
\begin{equation}
  \label{eq:2}
  \alpha=\frac23,\qquad\beta=\frac13,
\end{equation}
and with this choice (\ref{eq:1}) converges to 
\begin{equation}
  \label{eq:3}
  \frac{\partial}{\partial y}\psi=H\psi,\qquad 
  H=-\frac{\partial^2}{\partial u^2}+u
\end{equation}
regarded as a self-adjoint operator on $L^2(\IR)$. $H$ is the Airy
operator. The limit density field must correspond to free fermions
with $H$ as one-particle Hamiltonian. The fermions are in their ground
state at zero chemical potential.

Let us first describe the Fermi field in more detail. The Airy
operator $H$ has
$\IR$ as spectrum, which is purely absolutely continuous. The
generalized eigenfunctions are the Airy functions,
\begin{equation}\label{d.f}
- \frac{d^2}{du^2} \Ai (u - \lambda) + u \Ai(u- \lambda) = \lambda
\Ai (u - \lambda)\,.
\end{equation}
In particular the completeness relation
\begin{equation}\label{d.g}
\int d \lambda\, \Ai (u - \lambda) \Ai (v - \lambda) = \delta (u-v)\,
\end{equation}
holds. $K$ denotes the spectral projection onto $\{H\leq0\}$. Its
integral kernel is the Airy kernel
\begin{eqnarray}\label{d.h}
K(u, v) &=& \int^0_{- \infty} d \lambda\, \Ai (u-\lambda) \Ai (v -
\lambda) \nonumber \\
&=& \frac{1}{u-v} \big(\Ai (u) \Ai'(v) - \Ai'(u) \Ai(v)\big)\,.
\end{eqnarray}
Next we introduce the Fermi field  $a(u)$, $a^*(u)$, indexed
by $u\in\IR$. To distinguish from the fermions on a lattice we should
use a different symbol. Since the latter will not reappear, we find it
more convenient to stick to familiar notation. Integrated over a
test function $f\in L^2(\IR)$ the Fermi field becomes 
$a(f)=\int du\,f^*(u)a(u)$, $a^*(f)=\int du\,f(u)a^*(u)=a(f)^*$. They
satisfy the canonical anticommutation relations 
$\{a(f),a(g)\}=0=\{a^*(f),a^*(g)\}$ and $\{a(f),a^*(g)\}=(f,g)$ with
$(\cdot,\cdot)$ denoting the inner product of $f\in{L}^2(\IR)$
\cite{Brat}, and generate the CAR algebra ${\cal A}$. On ${\cal A}$ we
define the quasifree state 
$\omega$ through $\omega\big(a(f)\big)=0=\omega\big(a^*(f)\big)$ and
$\omega\big(a^*(f)a(g)\big)=(f,K\,g)$. In
particular, the moments of the density field are given by 
\begin{equation}
  \label{eq:9}
  \omega\bigg(\prod_{n=1}^ma^*(u_k)a(u_k)\bigg)=\det
  K(u_k,u_{k'})_{1\leq k,k'\leq m},
\end{equation}
for pairwise distinct $u_1,\dots,u_m$, compare with
(\ref{c.g}). (\ref{eq:9}) is the $m$-th correlation 
function. It vanishes at coinciding points.

To extend to unequal times one defines the Euclidean propagator 
\begin{eqnarray}
  \label{eq:9c}
  K(u,y;u',y')&=&\Big(e^{-yH}\big(K-\id\theta(y-y')\big)e^{y'H}\Big)
  (u,u')
  \nonumber\\
  &=&\text{sign}(y'-y)\int
  d\lambda\theta\big(\lambda(y-y')\big)e^{\lambda(y'-y)} 
  \Ai(u-\lambda)\Ai(v-\lambda),
\end{eqnarray}
for $y\neq y'$, written in terms of eigenfunctions of $H$, and
\begin{equation}
  \label{eq:9c1}
  K(u,y;u,y)=K(u,u),
\end{equation}
compare with (\ref{eq:prop2}), (\ref{eq:prop2a}). 
The quasifree state $\omega$ and the propagator are both determined by
the Airy operator $H$, which implies that $K$ depends only on $y-y'$.

The Airy field, denoted by $\xi(f,y)=\int du\, f(u)\xi(u,y)$, is the
density field of the Fermi   system defined through (\ref{eq:9c}),
(\ref{eq:9c1}). As in (\ref{c.r}), its moments are of determinantal
form and given by 
\begin{equation}
  \label{eq:9d}
  \IE\bigg(\prod_{k=1}^m\xi(f_k,y_k)\bigg)=
  \int\prod_{k=1}^mdu_k\,f_k(u_k)\,\,
  \det\big(K(u_k,y_k;u_{k'},y_{k'})\big)_{1\leq k,k'\leq m}.
\end{equation}
As it stands the left hand side of (\ref{eq:9d}) is only defined for
pairwise distinct 
$y_1,\dots,y_m$, but as in (\ref{c.r}) it can be continuously extended
to arbitrary time arguments.
Since $K$ depends only on
$y-y'$, the Airy field $\xi(f,y)$ is stationary in $y$.

Having introduced the limit object we turn to the edge scaling.
Recall that the PNG droplet has curvature. Therefore 
we set the scaled density field of the multi-layer PNG model, denoted
by $\xi_t$, as
\begin{equation}\label{d.b}
\xi_t (u,y) = t^{1/3}\eta_t \big([2t + t^{1/3}(u-y^2)], t^{2/3}y \big)\,,
\end{equation}
$[\,\cdot\,]$ denoting the integer part.
When integrated over the real, smooth, and rapidly decreasing test
function $f$ we have 
\begin{eqnarray}\label{d.c}
\xi_t(f,y) &=& \int du f(u) \xi_t(u,y) \nonumber \\
&=& \sum^{\infty}_{j = -\infty} f\big(t^{-1/3}
(j-2t)+y^2\big)\eta_t(j,t^{2/3}y ) + \Ord(t^{-1/3}) \nonumber\\
&=& \sum_{\ell\leq0} f
\big(t^{-1/3}(h_\ell(t^{2/3}y, t)-2t )+y^2 \big) +
\Ord(t^{-1/3})\,,
\end{eqnarray}
where the error of order $t^{-1/3}$ results from integrating over
cells of size $t^{-1/3}$ in the defining identity. Thus (\ref{d.c})
shows that through controlling the limiting moments of $\xi_t(f,y)$
one can infer the limit of the scaled height lines
$t^{-1/3}\big(h_\ell(t^{2/3}y,t)-2t\big)+y^2$. 

\begin{theorem}\label{the4}
Let $f_1, \dots, f_m$ be smooth test functions of compact support.
Then the following
limit holds,
\begin{equation}\label{d.p}
\lim_{t \rightarrow \infty} \IE_t \bigg(\prod^m_{k =1} \xi_t (f_k,
y_k) \bigg) = \IE \bigg(\prod^m_{k=1} \xi (f_k, y_k) \bigg)\,.
\end{equation}
\end{theorem}
{\it Proof:} 
Comparing (\ref{c.r}) and (\ref{eq:9d}) we have to establish that the
propagator (\ref{eq:prop2}), properly scaled, converges to the continuum
propagator (\ref{eq:9c}). This limit can be handled most directly in
the representation (\ref{eq:prop3}). We will need a separate argument for
$y\neq y'$ and for the left, resp. right, limit $y=y'$.

The case $y<y'$ runs in complete parallel to $y>y'$. To simplify
notation let us assume $y<y'$. The propagator for the scaled density
field is  
\begin{eqnarray}
  \label{d.p1}
  K_t(u,y;u',y')&=&e^{-2t^{2/3}y}e^{(u-y^2)y}
  t^{1/3}B_t\big([2t+t^{1/3}(u-y^2)],t^{2/3}y,
  [2t+t^{1/3}(u'-{y'}^2)],t^{2/3}y'\big) 
  \nonumber\\
  &&\times
  e^{2t^{2/3}y'}e^{-(u'-{y'}^2)y'}.
\end{eqnarray}
$y$, $y'$ are fixed and $K_t$ is
considered as a function on $\IR^2$. It is constant over cells of size
$t^{-1/3}$.
We used here the freedom that the determinant of (\ref{c.r}) does not
change under a similarity transformation and multiplied with the
factor
$\exp\big(2t^{2/3}(y'-y)\big)\exp\big((u-y^2)y-(u'-{y'}^2)y'\big)$
which diverges as $t\to\infty$. $K_t$, with the obvious
extension to  $y>y'$, determines the moments of the scaled density
field as
\begin{equation}
  \label{eq:d.p1a}
  \IE_t\bigg(\prod_{k=1}^m\xi_t(f_k,y_k)\bigg)=
  \int\prod_{k=1}^mdu_k\,f_k(u_k)\,\,
  \det\big(K_t(u_k,y_k;u_{k'},y_{k'})\big)_{1\leq k,k'\leq m}+\Ord(t^{-1/3}),
\end{equation}
clearly analoguous to (\ref{eq:9d}).

We insert (\ref{eq:prop3}) into (\ref{d.p1}). Then 
\begin{eqnarray}
  \label{d.p2}
  K_t(u,y;u',y')&=&t^{-1/3}\sum_{l\in t^{-1/3}\IN_-}e^{(y'-y)l}t^{1/3}
  J_{[2t+t^{1/3}(u-y^2-l)]}\big(2t\sqrt{1-t^{2/3}y^2}\big)
  \nonumber\\
  &&\times t^{1/3}J_{[2t+t^{1/3}(u'-{y'}^2-l)]}
  \big(2t\sqrt{1-t^{2/3}{y'}^2}\big)
  \nonumber\\
  &&\times \bigg\{\exp\big({-(2t^{2/3}y-(u-y^2)y-ly)}\big)
  \bigg(\frac{1+t^{-1/3}y}{1-t^{-1/3}y}\bigg)^{\left(2t+
      t^{1/3}(u-y^2-l)\right)/2}
  \nonumber\\
  &&\times \exp\big({(2t^{2/3}y'-(u-{y'}^2)y'-ly')}\big)
  \bigg(\frac{1-t^{-1/3}y'}{1+t^{-1/3}y'}\bigg)^{\left(2t+
      t^{1/3}(u'-{y'}^2-l)\right)/2}\bigg\}
  \nonumber\\
\end{eqnarray}
If $y<y'$ and $l\leq0$, the term $\big\{\cdots\big\}$ is uniformly
bounded in $t$, $l$ and converges to $1$ as $t\to\infty$. In fact this
holds uniformly in $u$, $u'$ on compact sets. By a result of Landau
\cite{Landau} 
\begin{equation}
  \label{d.p3}
  \sup_nt^{1/3}|J_n(2t)|\leq c/2^{1/3}
\end{equation}
with $c=0.7857\cdots$. From the asymptotics of integer Bessel
functions, cf.~(9.3.23) of \cite{Abra}, we conclude, uniformly for $u$
varying over a compact set,
\begin{equation}
  \label{d.p4}
  \lim_{t\to\infty}J_{[2t+t^{1/3}(u-y^2)]}\big(2t\sqrt{1-t^{2/3}y^2}\big)
  =\Ai(u).
\end{equation}
Since $e^{(y'-y)l}$, $l\leq0$, is integrable, by dominated convergence
\begin{equation}
  \label{d.p5}
  \lim_{t\to\infty}K_t(u,y;u',y')=K(u,y;u',y')
\end{equation}
uniformly over compact $u,u'$ sets.

Next we consider $y'\searrow y$ in (\ref{d.p1}), the right hand limit
being handled analogously. If the discrete Bessel kernel is
transformed according to (\ref{c.ad}), then
\begin{eqnarray}
  \label{eq:d.p6}
  K_t(u,y;v,y)&=&e^{y(u-u')}\bigg(\frac{1+t^{-1/3}y}
  {1-t^{-1/3}y}\bigg)^{[2t+t^{1/3}(u-y^2)]/2} \bigg(\frac{1-t^{-1/3}y}
  {1+t^{-1/3}y}\bigg)^{[2t+t^{1/3}(v-y^2)]/2} 
  \nonumber\\
  &&\times
  t^{1/3}B_{t\sqrt{1-t^{-2/3}y^2}}\big([2t+t^{1/3}(u-y^2)],
  [2t+t^{1/3}(v-y^2)]\big).
\end{eqnarray}
The first factor is uniformly bounded over compact $u,v$ sets and
converges to $1$ as $t\to\infty$. By Proposition (4.1) of \cite{BOO}
the discrete Bessel kernel with our scaling converges to the Airy
kernel $K(u,v)$ uniformly on compact $u,v$ sets as $t\to\infty$.

We conclude that (\ref{d.p5}) holds not only for $y\neq y'$ but also
for its right and left limits. Our claim follows by taking the limit
$t\to\infty$ in (\ref{eq:d.p1a}) which then yields
(\ref{eq:9d}). $\Box$\medskip 

The Airy field $\xi(f,y)$ is stationary in $y$. $\xi(f,y)$ is a point
process for fixed $y$. Its average density is given by
\begin{equation}\label{d.aj}
\IE(\xi(u,y)) = -u\, \Ai(u)^2 + \Ai' (u)^2
\end{equation}
which has the asymptotics \cite{For}
\begin{equation}\label{d.ak}
\IE\big(\xi(u,y)\big)\simeq\left\{
    \begin{array}{cl}
      \frac1\pi|u|^{1/2}-\frac1{4\pi|u|}\cos(4|u|^{3/2}/3)
      +\Ord\big(|u|^{-5/2}\big)& 
      \mbox{for $u\to-\infty$},\\[3mm]
      \frac{17}{96\pi}u^{-1/2}\exp\big(-4u^{3/2}/3\big)&\mbox{for
      $u\to\infty$}. 
    \end{array}
    \right.
\end{equation}
Note that for $u \rightarrow \infty$ the density decays quickly
because of the increasing linear potential, whereas for $u
\rightarrow - \infty$ the density is limited through the Fermi
exclusion. In particular, the point process for $\eta(f,y)$ has a
last point at $h_0(y)$ with probability one. Since all points are
distinct \cite{Sosh}, one can label as
\begin{equation}\label{d.al}
h_0(y) > h_{-1}(y) > \dots \,\,.
\end{equation}
$y \mapsto h_\ell(y)$, $\ell \in \IN_-$, are the fermionic world lines
underlying the Airy field. As to be shown in Appendix A $y
\mapsto h_\ell(y)$ is continuous with probability one. Moreover
$\langle\big(h_\ell(y)-h_\ell(y')\big)^2\rangle\simeq2|y-y'|$, which suggests
that the path measure for $\{h_\ell(y)$, $|y|<c,\,\ell=-M,\dots,0\}$
is absolutely continuous with respect to the Wiener measure,
i.e.~locally $h_\ell(y)$ is a modified Brownian motion.

Our main focus is the last fermion line.
\begin{definition}\label{def2}
Let $\xi(f,y)$ be the Airy field. The last world line, $h_0(y)$, is
called the Airy process and denoted by $A(y)$.
\end{definition}
We collect the basic properties of the Airy process.
\begin{theorem}\label{the5}
The Airy process $A(y)$ has continuous sample paths. $A(y)$ is
stationary. For given $y$, $A(y)$ has the distribution of $\chi_2$
of Tracy-Widom, see below eq.~(\ref{a.a}).
\end{theorem}

The convergence of the multi-layer PNG model to the Airy field
implies that the shape fluctuations of the PNG droplet converge to
the Airy process as $t 
\rightarrow \infty$. The following theorem is the precise version
of the main result, Theorem \ref{the0}, stated in the
Introduction. 
\begin{theorem}\label{the6}
Let $h(x,t)$ be the height of the PNG droplet and $h_t(y)$ its scaled
version according to (\ref{a.c}). Let $A(y)$ be the Airy process. Then
for any $m$, $y_j$,  $a_j \in \IR$, $j = 1, \dots, m$, we have
\begin{equation}\label{d.am}
  \lim_{t \rightarrow \infty} \IP_t \big(\{h_t (y_j) + y^2_j \leq
  a_j, j=1, \dots, m \}\big)=\IP \big( \{A(y_j) \leq a_j, j = 1, \dots,
  m \}\big)\,. 
\end{equation}
\end{theorem}
{\it Proof}: Let $f_j$ be the indicator of the interval $(a_j,
\infty)$. Then (\ref{d.am}) means
\begin{equation}\label{d.an}
\lim_{t \rightarrow \infty} \IP_t \Big( \bigcap^m_{j=1} \{\xi_t
(f_j, y_j) = 0 \}\Big) = \IP \Big( \bigcap^m_{j=1} \{\xi(f_j,
y_j) = 0 \} \Big)\,.
\end{equation}
We choose $a$ sufficiently large and split as $f_j = f^a_j + g^a$,
where $f^a_j$ is the indicator function of the interval $(a_j, a]$
and $g^a$ is the one of $(a, \infty)$. Then $\xi_t(f_j, y_j) =
\xi_t (f^a_j, y_j) + \xi_t(g^a, y_j)$. By Theorem \ref{the4} the joint
moments of $\xi_t (f^a_j, y_j)$, $j = 1, \dots, m$, converge to
their limit. Since their limit measure on $\IR^m$ is uniquely
defined by its moments we conclude that (\ref{d.an}) holds with
$f_j$ replaced by $f^a_j$. Up to constants the error term is
bounded by a sum over terms of the form 
\begin{equation}\label{d.ao}
  \IP_t \big(\xi_t (g^a, y_j) \geq 1 \big)\leq\IE_t \big(\xi_t
  (g^a, y_j)\big)=\sum_{j \geq t^{1/3}a} B_{t \sqrt{1 - t^{-2/3}y^2_j}} (j,j)
\end{equation}
which has a bound $C(a)$ uniform in $t$ such that $C(a) \rightarrow
0$ exponentially as $a \rightarrow \infty$, compare with
(\ref{d.ak}). $\Box$

\section{Some properties of the Airy process, two-point function}
\label{sec.e} 
\setcounter{equation}{0}

The scale invariant statistics of the PNG droplet is governed by
the Airy process. To gain some more quantitative information we
have to study the Airy process, most prominently its distribution
at a single point and at two points.

We denote by $P_a$ the projection onto the interval $(a, \infty)$,
i.e.~$P_a \psi(u) = \chi_a (u) \psi (u)$ with $\chi_a(u)$ the
indicator function of the set $(a, \infty)$. By definition
\begin{equation}\label{e.a}
\IP \big(\{A(y) \leq a \}\big) = \IP \big(\{\xi (\chi_a, y)
= 0 \}\big)\,.
\end{equation}
Let us set $\widehat{N}(\chi_a)=\int_a^\infty a^*(u)a(u)du$. Then 
\begin{equation}
  \label{eq:contprob}
  \IP(\xi(\chi_a,y)=0)=\lim_{\beta\to\infty}\omega\big(e^{-\beta
  \widehat{N}(\chi_a)}\big). 
\end{equation}
Since $\omega$ is quasifree,
\begin{equation}
  \label{eq:contdet}
  \omega\big(e^{-\beta \widehat{N}(\chi_a)}\big)=\det[1+(e^{-\beta}-1)P_aK],
\end{equation}
where the determinant is in $L^2(\IR)$. $P_aK$ is of trace class
\cite{Tra} and taking $\beta \rightarrow
\infty$ yields
\begin{equation}\label{e.d}
\IP \big(\{A(y) \leq a \}\big) = \det [1- P_aK ]\,.
\end{equation}
This determinant is studied in \cite{Tra} and shown to be related
to the Painlev\'e II differential equation. A plot for the
probability distribution of $A(y)$ can be found, for
  example, in \cite{Pra,Pra2}.

The next quantity of interest is the joint distribution of $A(0)$,
$A(y)$, where by reversibility it suffices to consider $y>0$. By the
same scheme as before one computes 
\begin{eqnarray}\label{e.e}
  \IP \big(\{A(0) \leq a, A(y) \leq b \}\big) &=& \IP \big(\{
  \xi(\chi_a, 0) = 0\}\cap\{\xi (\chi_b, y) = 0 \} \big)
  \nonumber \\
  &=&\lim_{\beta\to\infty}\omega\big(e^{-\beta\widehat{N}(\chi_a)}
  e^{-y\widehat{H}}e^{-\beta\widehat{N}(\chi_b)}e^{y\widehat{H}}\big)\,. 
\end{eqnarray}
Since $\omega$ is quasifree,
\begin{equation}
  \label{eq:conttwopt}
  \omega\big(e^{-\beta\widehat{N}(\chi_a)}e^{-y\widehat{H}}
  e^{-\beta\widehat{N}(\chi_b)}e^{y\widehat{H}}\big)=\det(1-B_\beta) 
\end{equation}
with
\begin{equation}
  \label{eq:bbeta}
  B_\beta=(1-e^{-\beta})(P_aK+e^{-yH}P_be^{yH}K)-(1-e^{-\beta})^2P_a
 e^{-yH}P_be^{yH}K.
\end{equation}
$e^{-yH}P_be^{yH}K$ is trace class, cf.~Appendix B. Thus
\begin{equation}\label{e.f}
\IP \big( \{A(0) \leq a, A(y) \leq b \}\big) = \det [1-B]
\end{equation}
and 
\begin{equation}
  \label{eq:b}
  B=P_aK+e^{-yH}P_be^{yH}K-P_ae^{-yH}P_be^{yH}K.
\end{equation}
Clearly, the determinant converges to 1 as $a, b \rightarrow
\infty$ and to 0 as $a, b \rightarrow - \infty$.

The properties of (\ref{e.d}) suggest that also the joint
distribution might satisfy a differential equation. We did not
succeed in finding one. Since the main interest is large $y$, we
rely on standard asymptotics by employing the expansion
\begin{equation}\label{e.g}
\log \det [1-B] = - \sum^{\infty}_{n=1} \frac{1}{n} \mathrm{tr}
[B^n]\,.
\end{equation}
When taking the trace of $B^n$, we see that $e^{yH}K $ is a
bounded operator, but $e^{-yH}$ remains unbalanced, since $H$ is
not bounded from below. E.g., $\mathrm{tr} [B]$ diverges as $\exp
(y^{3/2})$ for $y\to\infty$. The
form $\mathrm{tr}[B^n]$ is not suited for studying large $y$.

Such a situation is well known in the theory of Fermi systems
\cite{Matt, Salm}. Since the Dirac sea is filled up to energy
zero, one has to work with a new representation of the
CAR algebra, which
means to introduce field operators for the particles (energy $\geq
0$) and for the holes (energy $\leq 0$). In this representation
the Hamiltonian is positive. There is no need here to enter into
the full theory. It suffices to note that the series in
(\ref{e.g}) can be resummed such that only decaying exponentials
appear. Let $\sigma_n = \big(\sigma(1), \dots, \sigma(n)\big)$ be
an $n$-letter word where each letter is either $a$ or $b$,
$\sigma_0=\sigma_n$. Then
\begin{equation}\label{e.i}
\log \IP \big( \{ A(0) \leq a, A(y) \leq b \} \big) = -
\sum^{\infty}_{n =1} \frac{1}{n} \sum_{\{\sigma_n\}} \mathrm{tr}
\Big[\prod^n_{j=1} P_{\sigma(j)} K_{\sigma(j),\sigma(j+1)} \Big],
\end{equation}
$\sigma(n+1)\equiv\sigma(1)$,
with the following convention 
\begin{eqnarray}
  \label{eq:ei1}
  K_{a,a}&=&K_{b,b}=K
  \nonumber\\
  K_{a,b}&=&e^{-yH}(K-1)
  \nonumber\\
  K_{b,a}&=&e^{yH}K. 
\end{eqnarray}
The large $y$ asymptotics is extracted from (\ref{e.i}).

The spectral representation of $e^{-y H}(K-1) $, resp. $e^{yH}K$,
yields an integral of the form $- \int^{\infty}_0 d \lambda\,
e^{-\lambda y} g_+ (\lambda)$, resp. $\int^0_{- \infty} d \lambda\,
e^{\lambda y} g_- (\lambda)$, with some spectral functions $g_+$,
$g_-$. For large $y$ the weight concentrates at $\lambda = 0$ and
results in the asymptotics $- g(0)/y$, resp. $g(0)/y$. Therefore a
summand in (\ref{e.i}) decays as $y^{-\alpha}$, where $\alpha$ is
the number of broken bonds, i.e. $\dots ab \dots$ and $\dots ba
\dots$, in the word $\sigma_n$ of length $n$. The only words with
no broken bonds are either all $a$'s or all $b$'s. They sum to
\begin{equation}\label{e.j}
\det (1-P_aK ) \det(1- P_bK ).
\end{equation}
As to be expected, the Airy process is mixing and far apart events
become independent.

To order $y^{-2}$ we only allow two broken bonds which for large
$y$ leads to
\begin{eqnarray}\label{e.k}
y^{-2} \sum^{\infty}_{m=1} \sum^{\infty}_{n=1} \frac{1}{m+n}
\langle \Ai, (P_a K)^m P_a \Ai \rangle \langle \Ai, (P_b K)^n P_b \Ai
\rangle \nonumber \\
= y^{-2} \int^{\infty}_0 d \lambda\, \langle \Ai, P_a K (e^{\lambda} -
P_aK )^{-1} P_a \Ai \rangle \langle \Ai, P_b K  (e^{\lambda} - P_bK )^{-1}
P_b\Ai \rangle \,,
\end{eqnarray}
Here $\langle \Ai, \cdot\, \Ai \rangle$ means inner product in position
space with respect to the Airy function $\Ai(u)$, i.e. for some
operator $R$ with integral kernel $R(u,v)$,
$\langle\Ai,R\,\Ai\rangle=\int du\int dv\Ai(u)R(u,v)\Ai(v)$. To
compute the 
two-point function to leading order in $1/y$ we integrate the
probability measure with distribution function (\ref{e.f}) against
$a$ and $b$ and insert
to leading order from (\ref{e.j}), (\ref{e.k}). Using that $\det
(1-P_aK )$ has a good decay for $a \rightarrow - \infty$ and
$\langle \Ai, P_a K (e^{\lambda} - P_a K )^{-1}P_a \Ai \rangle$ has a good
decay for $a \rightarrow \infty$, one is allowed to integrate by parts
to obtain
\begin{eqnarray}\label{e.l}
\lefteqn{\langle A(0) A(y) \rangle - \langle A(0) \rangle \langle A(y)
\rangle=} 
\nonumber \\
&&\frac{1}{y^2} \int^{\infty}_0 d \lambda\, \Big[ \textstyle\int da
\det(1 -
P_aK) \langle \Ai, P_aK (e^{\lambda} - P_aK)^{-1} P_a \Ai \rangle
\Big]^2 + \Ord(y^{-4}).
\end{eqnarray}
The Airy process is positively correlated and has a slow decay as
$1/y^2$.

\section{Conclusions}
\setcounter{equation}{0}
\label{sec:conc}

The height statistics of the PNG droplet, for large $t$ and under the
scaling (\ref{a.c}), are given by the Airy
process. 
By universality such a result should be valid for any one-dimensional
growth model in the KPZ class. The only condition is that 
locally the macroscopic shape must have a
\emph{non-zero curvature}. If the interface is flat on the average,
other universal distributions will show up \cite{Pra}. In such a
situation, at present, no information on the multi-point statistics
is available. For example in the PNG droplet we could lift the
restriction that nucleation events are allowed only above the ground
layer $[-t,t]$. By translation invariance of the dynamics and the
initial condition, $h(x,0)=0$, we have $\langle h(x,t)\rangle=2t$ for
large $t$ and the process $y\mapsto
t^{-1/3}\big(h(yt^{2/3},t)-2t\big)=h^0_t(y)$ is stationary. For fixed
$y$, $h^0_t(y)$ converges to the GOE Tracy-Widom distribution
\cite{Pra2,BR}. The problem of the joint
distribution of $h^0_t(y_1)$, $h^0_t(y_2)$ remains open.

The Airy process contains a wealth of statistical information, which
cannot be resolved easily in Monte-Carlo simulations. In the standard
numerical experiment one merely considers the second moment of the
height differences. It is convenient to subtract the asymptotic mean
as $\bar h(x,t)=h(x,t)-2\sqrt{t^2-x^2}$. The quantity of interest is
then 
\begin{equation}
  \label{eq:c1}
  \langle\big(\bar h(x,t)-\bar h(0,t)\big)^2\rangle=G_t(x)
\end{equation}
for large $t$. Our main result says that $G_t(x)$ is of scaling form
and given by 
\begin{equation}
  \label{eq:c2}
  G_t(x)\simeq t^{2/3} g(t^{-2/3}x).
\end{equation}
The scaling function $g$ can be expressed through the two-point
function  of the Airy process as
\begin{equation}
  \label{eq:c3}
  g(y)=\langle\big(A(y)-A(0)\big)^2\rangle.
\end{equation}
For small $y$, one
has 
\begin{equation}
  \label{eq:c4}
  g(y)\simeq2|y|,\quad y\to0.
\end{equation}
On the other hand, for large $y$, $A(y)$ becomes independent from
$A(0)$ and
\begin{equation}
  \label{eq:c5}
  g(y)\simeq2a_2,\qquad|y|\to\infty
\end{equation}
with $a_2$ the truncated second moment of $\chi_2$,
$a_2=\langle\chi_2^2\rangle-\langle\chi_2\rangle^2\simeq0.81320$
numerically. The asymptotics (\ref{e.l}) says that 
\begin{equation}
  \label{eq:c6}
  g(y)\simeq2a_2-c|y|^{-2},\qquad |y|\to\infty,
\end{equation}
to next order. Because of the inverse operator in (\ref{e.l}) the
positive constant $c$ cannot be readily evaluated. A differential
equation, like the Painlev\'e II for the single point distribution,
would be helpful.
\begin{appendix}
  
\section*{Appendix A: Continuous sample paths of the Airy process}
\setcounter{equation}{0}
\label{appxa}
\renewcommand{\theequation}{A.\arabic{equation}}
\renewcommand{\thetheorem}{A.\arabic{theorem}}

Let $t\mapsto X(t)$ be a stochastic process with values in $\IR$. By a
criterion of Kolmogorov, cf.~\cite{Kal}, Theorem 2.23, if for some
constant $c$ 
\begin{equation}
  \label{eq:crit1}
  \IE(|X(t)-X(s)|^4)\leq c|t-s|^2,
\end{equation}
then $t\mapsto X(t)$ is continuous (in fact H\"older continuous with
exponent $<\frac14$) with probability one. Since (\ref{e.f}) provides
the joint distribution of $A(0)$, $A(y)$, (\ref{eq:crit1}) should be
an easy exercise. We did not succeed and rely on a more indirect
argument.

Let the test function $f$ be smooth and of compact support. By the
results of Section \ref{sec.d} we have 
\begin{equation}
  \label{eq:xi}
  \xi(f,y)=\sum_{j\leq0}f\big(h_\ell(y)\big).
\end{equation}
In the Lemma below we will prove that
\begin{equation}
  \label{eq:bound}
  \IE\big((\xi(f,y)-\xi(f,y'))^4\big)\leq c_f(y-y')^2.
\end{equation}
Therefore $y\mapsto\xi(f,y)$ is continuous with probability one. Since
the vague topology on locally finite point measures is countably
generated, the trajectory $y\mapsto\sum_{\ell\leq0}\delta(h_\ell(y)-u)$ is
continuous in the vague topology with probability one.
The convergence of a sequence of locally finite point measures in the
vague topology is equivalent to the convergence of each atom
\cite{Kerst}. Thus (\ref{eq:bound}) implies that $y\mapsto h_\ell(y)$
for each $\ell$ is
continuous with probability one. In particular
$A(y)=h_0(y)$ is continuous.
\begin{lemma}
  Let $f$ be smooth and of compact support. Then there is a constant
  $c_f$ such that
\begin{equation}
  \label{eq:bound2}
  \IE\big((\xi(f,y)-\xi(f,0))^4\big)\leq c_fy^2.
\end{equation}
\end{lemma}
{\it Proof:} As a warm-up,  and to fix notation, we first compute the
second moment. We suppress $f$ and set $\xi(f,y)=\xi_y$. We take
$y\geq0$. $y\leq0$ follows from reversal symmetry $y\mapsto-y$ and
stationarity in $y$. As shorthand we define $L=K-1$. $-L$ is the
projection operator onto ${H\geq0}$. We regard $f$ as a multiplication
operator, $f\psi(u)=f(u)\psi(u)$, and formally set
$f_y=e^{yH}f\,e^{-yH}$. By construction, in the expressions below,
only the bounded operators $Ke^{yH}$ and $e^{-yH}L$ appear as
factors. By the determinantal formula (\ref{eq:9d})
\begin{eqnarray}
  \label{eq:2cycle}
  \IE\big((\xi_y-\xi_0)^2\big)&=&2\big(\IE(\xi_0^2)-\IE(\xi_0\xi_y)\big)
  \nonumber\\
  &=&-2\big(\tr(fKfL)-\tr(fKf_yL)\big)=-2\tr(fK(f-f_y)L).
\end{eqnarray}
From the spectral representation, one concludes differentiability of
any order in $y$, $y>0$, and up to linear order
\begin{eqnarray}
  \label{eq:lo}
   \IE\big((\xi_y-\xi_0)^2\big)&=&2y\,\tr(fK[H,f]L)+\Ord(y^2)
   \nonumber\\
   &=&2y\big(\tr(fK[H,f]K)-\tr(fK[H,f])\big)+\Ord(y^2). 
\end{eqnarray}
For real operators one has $\tr(AB)=\big(\tr (AB)\big)^*=\tr(B^*A^*)$. Since
$[H,f]^*=-[H,f]$, the first summand vanishes and 
\begin{eqnarray}
  \label{eq:lo2}
   \IE\big((\xi_y-\xi_0)^2\big)&=&y\,\tr(fK[H,f]L)+\Ord(y^2)
   =y\,\tr(K[f,[H,f]])+\Ord(y^2)
   \nonumber\\&=&   2y\,\tr K({f'})^2+\Ord(y^2).
\end{eqnarray}
The variance (\ref{eq:lo2}) implies that
$\IE\big((h_\ell(y)-h_\ell(0))^2\big)=2|y|+\Ord(y^2)$ for each $\ell$ in the
limit $y\to0$, as to be expected from the construction of the Airy
field. 

The fourth moment requires more effort. We have, using stationarity
and reversibility
in $y$,
\begin{equation}
  \label{eq:fourth}
  \IE\big((\xi_y-\xi_0)^4\big)=2\big(\IE(\xi_0^4)
  -4\IE(\xi_0^3\xi_y)+3\IE(\xi_0^2\xi_y^2)\big).
\end{equation}
We again use the determinant formula (\ref{eq:9d}), which most
conveniently is decomposed into cycles.
For the fourth moment there are $4!$ permutations. They subdivide into
(i) four $1$-cycles (1 term), (ii) two $1$-cycles plus one
$2$-cycle (6 terms), (iii) two $2$-cycles (3 terms), (iv) one
$1$-cycle plus one $3$-cycle (8 terms), and (v) one $4$-cycle (6
terms). The sign of the permutation will be of no significance for the
argument.

(i) and (ii) vanish, (iii) is the ``Gaussian'' term,
\begin{equation}
  \label{eq:iii}
  \text{(iii)}=3\IE\big((\xi_y-\xi_0)^2\big)^2\leq c_f y^2
\end{equation}
from (\ref{eq:lo2}).
For (iv) we have a $1$-cycle and two $3$-cycles in reverse
order. Summing over all such cycles yields 
\begin{eqnarray}
  \label{eq:iv}
  \text{(iv)}&=&12\tr(Kf)\big(-\tr(fKf_yLfL)+\tr(fKf_yLf_yL)
  \nonumber\\
  &&\hspace{2cm}-\tr(fKfKf_yL)+\tr(fKf_yKf_yL)\big).
\end{eqnarray}
Since $\tr(fKf_yLfL)^*=\tr fKf_yLf_yL$ and
$\tr(fKfKf_yL)^*=\tr fKf_yKf_yL$ the contribution (iv) vanishes.

It remains to study the $4$-cycles. For their sum one obtains
\begin{eqnarray}
  \label{eq:v}
  \text{(v)}&=&\hspace{3mm}\tr(fKfKfK(f-f_y)L)-3\tr(fKfK(f-f_y)KfL)
  \nonumber\\
  &&+\tr(fKfK(f-f_y)LfL)-3\tr(fKfKfL(f-f_y)L)
  \nonumber\\
  &&+\tr(fKfK(f-f_y)LfL)-3\tr(fK(f-f_y)Kf_yLfL)
  \nonumber\\
  &&+\tr(fKfLfK(f-f_y)L)-3\tr(fK(f-f_y)LfKf_yL)
  \nonumber\\
  &&+\tr(fK(f-f_y)LfKfL)-3\tr(fKf_yLfK(f-f_y)L)
  \nonumber\\
  &&+\tr(fK(f-f_y)LfLfL)-3\tr(fKf_yL(f-f_y)LfL),
\end{eqnarray}
where we combined the terms such that the order $y$ is manifest. By
the spectral theorem (\ref{eq:v}) is differentiable to any order in
$y$, $y>0$. Thus to complete the proof it suffices to show that the
linear order of (\ref{eq:v}) vanishes. Since
$f-f_y\simeq-y[H,f]$, we obtain
\begin{eqnarray}
  \label{eq:lo3}
  \text{(v)}&=&y\,\tr\big([H,f](LfKfKfK-3KfLfKfK+2LfLfKfK
  -3LfKfKfL
  \nonumber\\
  &&-3KfLfLfK-4LfKfLfK+LfLfLfK-3LfLfKfL)\big)+\Ord(y^2).
\end{eqnarray}
We substitute $L=K-1$ and use that $[H,f]^*=-[H,f]$. The term with
four $K$'s reads
\begin{equation}
  \label{eq:fourK}
  -12\tr([H,f]KfKfKfK)=0.
\end{equation}
The term with three $K$'s reads
\begin{equation}
  \label{eq:threeK}
  -6\tr\big([H,f]fKfKfK+Kf^2KfK+KfKf^2K+KfKfKf)\big)=0,
\end{equation}
since first and fourth summand and second and third summand
cancel. The term with two $K$'s reads
\begin{equation}
  \label{eq:twoK}
  -\tr\big([H,f](6fKfKf+2Kf^3K+3fKf^2K+3Kf^2Kf)\big)=0,
\end{equation}
since the first and second summand vanish and since the third summand
cancels against the fourth one. 
The term with one $K$ reads, upon adding the adjoint,
\begin{equation}
  \label{eq:oneK}
  \tr([H,f]f^3-f^3[H,f]-3f[H,f]f^2+3f^2[H,f]f).
\end{equation}
At this point the specific structure of $H$ enters. We have
$[H,f]=-f''-2f'(d/du)$. Working out the commutators ensures also that
the last term vanishes. Thus the linear order from the $4$-cycles
vanishes. The quadratic order does not vanish, however, and reflects
the  deviations from the local Brownian motion statistics. $\Box$
\section*{Appendix B: Trace class properties}
\setcounter{equation}{0}
\label{appxb}
\renewcommand{\theequation}{B.\arabic{equation}}
\renewcommand{\thetheorem}{B.\arabic{theorem}}
As shown in \cite{Tra} $P_aK$ is of trace class. If one establishes
that $R=e^{-yH}P_be^{yH}K$ is of trace class, then each summand in
(\ref{eq:b}) is separately of trace class. In
the energy representation $R$ has the kernel
\begin{equation}
  \label{eq:Rkernel}
  R(\lambda',\lambda)=e^{-y\lambda'}\int_b^\infty
  dx\,\Ai(x-\lambda')\Ai(x-\lambda)e^{\lambda y}\chi(\lambda\leq0).
\end{equation}
If $\lambda\leq0$, $\lambda'\geq0$, the integral is bounded uniformly
in $\lambda,\lambda'$. If $\lambda\leq0$, $\lambda'\leq0$, one uses
that $\Ai(x-\lambda')\leq c\,e^{-\frac23|x-\lambda'|^{3/2}}$ for large
negative $\lambda'$. Thus the integral dominates the factor
$e^{-\lambda'y}$. This implies the bound 
\begin{equation}
  \label{eq:Rbound}
  |R(\lambda',\lambda)|\leq c\,e^{-\gamma(|\lambda'|+|\lambda|)}
\end{equation}
with suitable constants $c,\gamma>0$. Hence $R$ is of trace class.
\end{appendix}
\vspace{3mm}\\
{\bf Note added in proof:} In their recent preprint \cite{Okou} Okounkov
and Reshetikhin consider the $(1,1,1)$ interface of the
three-dimensional Ising model at zero temperature, which maps onto a
domino tiling of a form similar to the Aztec diamond. They prove that
correlations are of determinantal form and compute the limit
shape. Their Theorem~1 is the analogue of our Eq.~(\ref{c.r}).

\end{document}